\RequirePackage{fix-cm}
\documentclass[numbook,runningheads]{svjour3}
\usepackage{amssymb}
\usepackage{amsmath}
\usepackage{algorithm}
\usepackage{algpseudocode}
\usepackage{color}
\usepackage{graphicx}
\usepackage{hyperref}

\newenvironment{changes}{\color[rgb]{0.7,0,0}}{}
\providecommand{\mychanges}[1]{{\color[rgb]{0.7,0,0}#1}}
\renewenvironment{changes}{}{}
\renewcommand{\mychanges}[1]{{#1}}

\providecommand{\E}[1]{{\ensuremath{\mathrm{E}}\mspace{-2mu}\left[#1\right]}}
\providecommand{\prob}[1]{{\ensuremath{\mathrm{P}}\mspace{-2mu}\left[#1\right]}}
\providecommand{\var}[1]{{\ensuremath{\mathrm{Var}}\mspace{-2mu}\left[#1\right]}}

\providecommand{\tol}{\mathrm{TOL}}
\providecommand{\rset}{\mathbb{R}}
\providecommand{\zset}{\mathbb{Z}}
\providecommand{\nset}{\mathbb{N}}
\providecommand{\work}{\ensuremath{W}}

\providecommand{\est}[1]{\overset{\sim}{#1}}
\providecommand{\Order}[1]{ {\ensuremath{ \mathcal O\left( #1 \right)}} }
\providecommand{\order}[1]{ {\ensuremath{ o\left( #1 \right)}} }

\title{A Continuation {M}ultilevel {M}onte {C}arlo algorithm}

\author{Nathan~Collier
  \and Abdul--Lateef~Haji--Ali
  \and Fabio~Nobile
  \and Erik~von~Schwerin
  \and Ra\'{u}l~Tempone
}

\institute{
  N.~Collier  (\email{nathaniel.collier@gmail.com})\at {Oak Ridge National Lab, Climate Change Science Institute (CCSI), Environmental Sciences Division}
  \and
  A.~Haji--Ali (\email{abdullateef.hajiali@kaust.edu.sa}) \and R.~Tempone  (\email{raul.tempone@kaust.edu.sa})\at
        Applied Mathematics and Computational Sciences, KAUST, Thuwal, Saudi Arabia.
  \and
  F.~Nobile \and E.~Schwerin \at {MATHICSE-CSQI, EPF de Lausanne, Switzerland.} %
}

\def\makeheadbox{{}}
\date{}

\begin{document}

\maketitle

\begin{abstract}
  We propose a novel Continuation Multi Level Monte Carlo (CMLMC)
  algorithm for weak approximation of stochastic models. The CMLMC
  algorithm solves the given approximation problem for a sequence of
  decreasing tolerances, \mychanges{ending when the required error
    tolerance is satisfied.} CMLMC assumes discretization hierarchies
  that are defined a priori for each level and are geometrically
  refined across levels. The actual choice of computational work
  across levels is based on parametric models for the average cost per
  sample and the corresponding weak and strong errors. These
  parameters are calibrated using Bayesian estimation, taking
  particular notice of the deepest levels of the discretization
  hierarchy, where only few realizations are available to produce the
  estimates. The resulting CMLMC estimator exhibits a non-trivial
  splitting between bias and statistical contributions. We also show
  the asymptotic normality of the statistical error in the MLMC
  estimator and justify in this way our error estimate that allows
  prescribing both required accuracy and confidence in the final
  result. Numerical results substantiate the above results and
  illustrate the corresponding computational savings \mychanges{in
    examples that are described in terms of differential equations
    either driven by random measures or with random coefficients}.

\end{abstract}

\section{Introduction}\label{s:intro} % intro.tex
Multilevel Monte Carlo Sampling was first introduced for applications
in the context of parametric integration by
Heinrich~\cite{heinrich98,hs99}. Later, to consider weak approximation
of stochastic differential equations (SDEs) in mathematical finance,
Kebaier~\cite{kebaier05} introduced a two-level Monte Carlo technique
in which a coarse grid numerical approximation of an SDE was used as a
control variate to a fine grid numerical approximation, thus reducing
the number of samples needed on the fine grid and decreasing the total
computational burden. This idea was extended to a multilevel Monte
Carlo (MLMC) method by Giles in~\cite{giles08}, who introduced a full
hierarchy of discretizations with geometrically decreasing grid sizes.
By optimally choosing the number of samples on each level this MLMC
method decreases the computational burden, not only by a constant
factor as standard control variate techniques do, but even reducing
the rate in the computational complexity \begin{changes}to compute a
  solution with error tolerance $\tol > 0$\end{changes} from
$\Order{\tol^{-3}}$ of the standard Euler-Maruyama Monte Carlo method
to $\Order{\log{(\tol)}^2\tol^{-2}}$, assuming that the work to
generate a single realization is $\Order{\tol^{-1}}$. For
one-dimensional SDEs, the computational complexity of MLMC was further
reduced to $\Order{\tol^{-2}}$ by using the Milstein Scheme
\cite{giles08_m}. Moreover, the same computational complexity can be
achieved by using antithetic control variates with MLMC in
multi-dimensional SDEs with smooth and piecewise smooth payoffs
\cite{gs13}.

This standard MLMC method has since then been extended and applied in a wide variety of
contexts, including jump diffusions \cite{xg12} and Partial Differential Equations (PDEs) with random coefficients
\cite{bsz11,cst13,cgst11,gr12,tsgu13}.
It is \mychanges{shown in \cite[Theorem 2.3]{tsgu13}} that there is an optimal convergence rate that is similar to the previously mentioned complexity rates, but that
depends on the relation between
the rate of strong convergence of the discretization method of the
\mychanges{underlying equation} and the work complexity associated
with generating a single sample of the quantity of interest. In fact, in certain cases, the computational complexity can be of the
optimal rate, namely $\Order{\tol^{-2}}$.

To achieve the optimal MLMC complexity rate and to obtain an estimate of the statistical error,
   sufficiently accurate estimates of the variance on each level must be obtained.
Moreover, finding the optimal number of levels requires a sufficiently accurate estimate of the bias.
As such, an algorithm is needed to find these estimates without incurring a significant overhead to the estimation
of the wanted quantity of interest.
In \cite{giles08}, Giles proposed an algorithm, henceforth referred to as {\em Standard} MLMC or SMLMC,
   that works by iteratively increasing the number of levels and using sample variance estimates across levels.
Moreover, SMLMC uses an arbitrary fixed accuracy splitting between the bias and the statistical error contributions.
Other works \cite{teckentrup13,gs13c,gs13b,cgst11} listed similar versions of this algorithm.
We outline this algorithm in Section \ref{s:std}.

In Section \ref{s:adapt}, we propose a novel continuation type of MLMC
algorithm that uses models for strong and weak convergence and for
average computational work per sample. We refer to this algorithm as
Continuation MLMC or CMLMC. The CMLMC algorithm solves the given
problem for a sequence of decreasing tolerances, which plays the role
of a continuation parameter, \mychanges{the algorithm ends when the
  required error tolerance is satisfied}. Solving this sequence of
problems allows CMLMC to find increasingly accurate estimates of the
bias and variances on each level, in addition to the quantity of
interest, which is the goal of the computation. In each case,
\mychanges{given the current estimate of problem parameters, the
  optimal number of levels of the MLMC hierarchy is found}. Moreover,
we use a Bayesian inference approach to robustly estimate the various
problem parameters. The CMLMC algorithm is able to relax the
statistical error bound given the bias estimate, to achieve the
optimal splitting between the two. These techniques improve the
computational complexity of the CMLMC algorithm and decreases the
variability of the running time of the algorithm.

The outline of this work is as follows: We start in Section \ref{s:mlmc} by recalling the MLMC method and the assumed models on work, and on weak and strong convergence.
After introducing the algorithms in Sections \ref{s:std} and \ref{s:adapt}, Section \ref{s:res}
 presents numerical examples, which include three-dimensional PDEs
 with random inputs and an It\^o SDE.
Finally, we finish by offering conclusions and suggesting directions for future work in Section \ref{s:conc}.

\section{Multilevel Monte Carlo}\label{s:mlmc} % mlmc.tex
\subsection{Problem Setting}\label{sec:hier_intro}
Let $g(u)$ denote a real valued functional of the solution, $u$, of an underlying stochastic model.
We assume that $g$ is either a bounded linear functional or Lipschitz with respect to $u$.  Our goal is to
approximate the expected value, $\E{g(u)}$, to a given accuracy $\tol$
and a given confidence level. We assume that individual outcomes
of the
underlying solution $u$ and the evaluation of the functional $g(u)$ are approximated
by a discretization-based numerical scheme characterized by a mesh
size, $h$. The value of $h$ will govern the weak and strong errors in the
approximation of $g(u)$ as we will see below. To motivate this setting,
we now give two examples and identify  the numerical discretizations,
the discretization parameter, $h$, and the corresponding rates of approximation.
The first example is  common  in engineering applications like heat conduction and
groundwater flow.
The second example is a simple one-dimensional geometric Brownian motion with European call option.
\begin{example}
\label{ex:spde_problem}
Let $(\Omega,\mathcal{F},\mathbb P)$ be a complete probability space and $\mathcal{D}$ be a bounded convex polygonal domain in $\rset^d$.
 Find $u: \mathcal{D} \times \Omega \to \rset$ that solves almost surely (a.s.)
the following equation:
\begin{subequations}
\begin{align*}
    -\nabla \cdot \left( a(\mathbf x; \omega) \nabla u(\mathbf x; \omega) \right) &= f(\mathbf x; \omega) &\text{ for }
            \mathbf x  \in \mathcal D,  \\
     u(\mathbf x; \omega) &=  0 &\text{ for } \mathbf x \in \partial \mathcal D,
\end{align*}
\end{subequations}
\begin{changes}
  where $\omega \in \Omega$ and the value of the diffusion coefficient
  and the forcing are represented by random fields, yielding a random
  solution. We wish to compute $\E{g(u)}$ for some deterministic
  functional $g$ which is globally Lipschitz satisfying $|g(u) - g(v)|
  \leq G \|u-v\|_{H^1(\mathcal D)}$ for some constant $G>0$ and all
  $u,v\in H^1(\mathcal D)$.
Following \cite{tsgu13}, we also make the following assumptions
\begin{itemize}
  \item $a_{\min}(\omega) = \min_{\mathbf x \in \mathcal  D}a(\mathbf x; \omega) > 0$ a.s. and $1/a_{\min} \in  L^p_\mathbb{P}(\Omega)$, for all $p \in (0, \infty)$.
  \item $a \in L^p_{\mathbb P}(\Omega, C^1(\overline{\mathcal D}))$,
    for all $p \in (0, \infty)$.
  \item $f \in L^{p^*}_{\mathbb P}(\Omega, L^2(\mathcal D))$ for some
    $p^* > 2$.
  \end{itemize}
  Here, $L^p_{\mathbb P}(\Omega, \mathcal B)$ is the space of $\mathcal
  B$-valued random fields with a finite $p$'th moment of their
  $\mathcal B$-norm, where the $p$-moment is with respect to
  measure $\mathbb P$. On the other hand, $C^1(\overline D)$ is the
  space of continuously differentiable functions with the usual norm
  \cite{cst13}.
  Note that with these assumptions and since $\mathcal D$ is bounded,
  one can show that   $\max_{\mathbf x \in \mathcal  D}a(\mathbf x;
  \omega) < \infty$ a.s.
A standard approach to approximate the solution of the previous problem is to use Finite Elements on regular triangulation. In such a setting,
the parameter $h>0$ refers to either the maximum element diameter or another characteristic length and the corresponding
approximate solution is  denoted by $u_h(\omega)$.
For piecewise linear or piecewise $d$-multilinear continuous finite element
approximations, and with the previous assumptions,
it can be shown \cite[Corollary 3.1]{tsgu13} that asymptotically as
$h \to 0$:
\begin{itemize}
\item $\left|\E{g(u)-g(u_h)} \right| \lesssim Q_W\, h^2$ for a constant $Q_W>0$.
\item $\var{g(u)-g(u_h)} \lesssim Q_S\, h^4$ for a constant $Q_S > 0$.
\end{itemize}
\end{changes}

\end{example}

\begin{example}
\label{ex:sde_problem}
 Here we study the weak
approximation of It\^o stochastic differential equations (SDEs),
\begin{equation}
  du(t) = a(t,u(t))dt +
  b(t,u(t))dW(t), \qquad 0<t<T,
  \label{eq:sde}
\end{equation}
where $u(t;\omega)$ is a stochastic process in $\rset^d$, with
randomness generated by a $k$-dimensional Wiener process with
independent components, $W(t;\omega)$, cf.~\cite{KS,Ok}, and
$a(t,u)\in\rset^d$ and $b(t,u)\in\rset^{d\times k}$ are the drift and
diffusion fluxes, respectively.  For any given sufficiently well behaved function,
$g: \rset^d\rightarrow \rset$, our goal is to approximate the expected
value, $\E{g(u(T))}$. A typical application is to compute option prices in
mathematical finance, cf.~\cite{JouCviMus_Book,Glasserman}, and other
related models based on stochastic dynamics.
When  one uses a standard Euler Maruyama (Forward Euler) method based on  uniform
time steps of size $h$ to approximate \eqref{eq:sde},
the following rates of approximation hold: $|\E{g(u(T))-g(u_h(T))}| = Q_W \,h + \order{h}$ and
$\E{(g(u(T))-g(u_h(T)))^2} = Q_S \,h + \order{h},$ for some constants,
$0<Q_W,Q_S<\infty$\begin{changes}, different from the constants of the previous example.\end{changes}
For suitable assumptions on
the functions $a$, $b$ and $g$, we refer to \cite{Mordeckietal08,stz2001}.
\end{example}

To avoid cluttering the notation, we  omit the reference to the underlying solution from
now on, simply denoting the quantity of interest by $g$. Following the standard MLMC approach, we assume,  for any given non-negative integer $L \in \nset$, that we have a hierarchy of $L+1$
meshes defined by a decreasing sequence of mesh sizes $\{h_\ell\}_{\ell=0}^L$
where ${h_\ell=h_0 \beta^{-\ell}}$ for some $h_0>0$ and a constant integer $\beta > 1$.
We denote the resulting approximation of $g$ using mesh size $h_\ell$ by
$g_\ell$, or by $g_\ell(\omega)$ when we want to stress the dependence
on an outcome of the underlying random model.
Using the following notation:
\begin{align*}
    G_\ell(\omega) &=
    \begin{cases}
        g_0(\omega) & \text{if } \ell=0, \\
        g_\ell(\omega) - g_{\ell-1}(\omega) & \text{if } \ell >0,
    \end{cases}
  \end{align*}
the expected value of the finest approximation, $g_L$, can be expressed as
\begin{align*}
  \E{g_L} & = \sum_{\ell=0}^L \E{G_{\ell}},
\end{align*}
where the MLMC estimator is obtained by replacing the expected values
in the telescoping sum by sample averages.
\begin{changes}
We denote the sample averages by $\est{G}_\ell$ as
\begin{align*}
  \est G_\ell & =  M_\ell^{-1}\sum_{m=1}^{M_\ell} G_\ell(\omega_{\ell,m}).
\end{align*}
\end{changes}
\mychanges{Each sample average, $\est G_\ell$, is computed using
  $M_\ell \in \zset_+$ independent identically
distributed (i.i.d.) outcomes, $\{\omega_{\ell,m}\}_{m=1}^{M_\ell}$, of
the underlying, mesh-independent, stochastic model; i.e.
$\omega_{\ell,m} \in \Omega$ for all $\ell$ and $m$}.
The MLMC estimator can then be written as
\begin{align}
  \label{eq:MLMC_estimator}
  \mathcal{A} & = \sum_{\ell=0}^L \est G_\ell.
\end{align}
Note that the outcomes are also assumed to be independent among the
different sample averages, $\{\est G_\ell\}_{\ell=0}^L$.

We use the following model for the expected value of the cost
associated with generating one sample of $G_{\ell}$, including
generating all the underlying random variables:
\begin{changes}
\begin{equation*}
W_\ell \propto h_\ell^{-\gamma} = h_0^{-\gamma}\beta^{\ell \gamma}
\end{equation*}
\end{changes}
for a given $\gamma$. Note the cost of generating a sample of $G_\ell$ might differ for different realizations, for example due to different number of iterations in an
iterative method or due to adaptivity of the used numerical method.
The parameter $\gamma$ depends on the number of dimensions of the underlying problem and the used numerical method.
For example, $\gamma=1$ for the one-dimensional SDE in Example \ref{ex:sde_problem}.
For the PDE in Example \ref{ex:spde_problem}, if the number of dimensions is $d=3$ then $\gamma=3 \tilde \gamma$, where $\tilde \gamma$ depends on the solver used to solve the resulting linear system.
In that example, iterative methods may have a smaller value of $\tilde \gamma$ than direct methods.
The theoretical best-case scenario for iterative methods would be $\tilde \gamma = 1$ for multigrid methods.
On the other hand, we would have $\tilde \gamma = 3$  if one used a direct method using a naive Gaussian elimination on dense matrices.
The total work of the estimator~\eqref{eq:MLMC_estimator} is
\begin{align*}
  \work & = \sum_{\ell=0}^{L} M_\ell W_\ell.
\end{align*}

We want our estimator to satisfy a tolerance with prescribed failure probability $0 < \alpha < 1$, i.e.,
\begin{align}
  \label{eq:goal}
  \prob{\left|\E{g}-\mathcal{A}\right| > \tol} &\leq \alpha,
\end{align}
while minimizing the work, $\work$.  Here, we split the total error into bias and statistical error,
\begin{align*}
    \left|\E{g}-\mathcal{A}\right| \leq \underbrace{\left|\E{g-\mathcal{A}}\right|}_{\text{Bias}}
    + \underbrace{\left|\E{\mathcal{A}}-\mathcal{A}\right|}_{\text{Statistical error}},
\end{align*}
and use a splitting parameter, $\theta \in (0,1)$, such that
\begin{changes}
\[ \tol =  \underbrace{(1-\theta) \tol}_{\text{Bias tolerance}} +
\underbrace{\hskip 1cm \theta \tol \hskip1cm }_{\text{Statistical error tolerance}}.\]
\end{changes}
The MLMC algorithm should bound the bias, $B = \left|\E{g-\mathcal{A}}\right|$, and the statistical error as follows:
\begin{subequations}
\begin{align}
    B =  \left|\E{g-\mathcal{A}}\right| &\leq (1-\theta)\tol,  \\
  \left|\E{\mathcal{A}}-\mathcal{A}\right| &\leq \theta \tol,
  \label{eq:stat_error}
\end{align}
\end{subequations}
where the latter bound should hold with probability $1-\alpha$. Note that $\theta$ does not have to be a constant, indeed it can depend on $\tol$ as we shall
see in Section \ref{s:adapt}.
In the literature, some authors (e.g. \cite{giles08}) have controlled the mean square error (MSE),
    \[\text{MSE} = {\left|\E{g-\mathcal{A}}\right|}^2 + \E{\left|\E{\mathcal{A}}-\mathcal{A}\right|^2}, \]
 rather than working with \eqref{eq:goal}. We prefer to work with \eqref{eq:goal} since it allows us to prescribe both the accuracy $\tol$ and the
 confidence level, $1-\alpha$, in our results.
The bound \eqref{eq:stat_error} leads us to require
\begin{align}
  \label{eq:var_bound}
  \var{\mathcal{A}} & \leq \left(\frac{\theta \tol}{C_\alpha} \right)^2,
\end{align}
for some given confidence parameter, $C_\alpha$, such that ${\Phi(C_\alpha) = 1-\frac{\alpha}{2}}$; here, $\Phi$ is the cumulative distribution
function of a standard normal random variable.
The bound \eqref{eq:var_bound} is motivated by the Lindeberg Central Limit Theorem in the limit ${\tol\to 0}$, cf. Lemma \ref{thm:clt_result} in the Appendix.

By construction of the MLMC estimator, $\E{\mathcal{A}}=\E{g_L}$, and denoting $ V_\ell = \var{G_\ell}$, then by independence, we have
  $\var{\mathcal{A}} = \sum_{\ell=0}^L V_\ell M_\ell^{-1},$
and the total error estimate can be written as
\begin{align}
    \text{Total error estimate} = B + C_\alpha \sqrt{\var{\mathcal{A}}}.
    \label{eq:total-error}
\end{align}
Given $L$ and $0 < \theta < 1$ and minimizing $W$ subject to the statistical constraint \eqref{eq:var_bound} for
        $\left\{M_\ell\right\}_{\ell=0}^L \in \rset^{L+1}$ gives the following optimal number of samples per level $\ell$:
\begin{equation}
    M_\ell = \left( \frac{C_\alpha}{\theta \tol} \right)^2 \sqrt{\frac{V_\ell}{W_\ell}} \left( \sum_{\ell=0}^L \sqrt{V_\ell W_\ell} \right).
    \label{eq:optimal_ml}
\end{equation}
When substituting the optimal number of samples in all levels the optimal work can be written in terms of $L$ as follows
\begin{equation}
    \work(\tol, L) = \left( \frac{C_\alpha}{\theta \tol} \right)^2 \left(\sum_{\ell=0}^L \sqrt{V_\ell W_\ell}\right)^2.
    \label{eq:work_per_L}
\end{equation}
Of course, the number of samples on each level is a positive integer.
To obtain an approximate value of the optimal integer number of samples, we take the ceiling of the real-valued optimal values in \eqref{eq:optimal_ml}.

In this work, we assume the following models on the weak error and variance:
\begin{subequations}
\begin{align}
    \label{eq:weak_error_model} \E{g-g_\ell} &\approx Q_W h_\ell^{q_1},\\
    \var{g_\ell-g_{\ell-1}} &\approx Q_S h_{\ell-1}^{q_2},
\end{align}
\end{subequations}
for some constants $Q_W\neq 0,Q_S>0,q_1>0$ and $0 < q_2 \leq 2 q_1$.
For example, recall that the PDE in Example \ref{ex:spde_problem} has $q_2=2 q_1$
and in Section \ref{s:res}, the PDE is solved using a finite element method with standard
trilinear basis and it has $q_1=2$.
On the other hand, for the SDE in Example \ref{ex:sde_problem} with Euler discretization, $q_1=q_2=1$.
Collectively, we refer to the parameters $q_1, q_2, Q_S, Q_W$ and $\{V_\ell\}_{\ell=0}^L$ as problem parameters.
Based on these models, we can write for $\ell>0$
\begin{changes}
\begin{subequations}
\begin{align}
    \E{G_\ell} &\approx Q_W h_0^{q_1} \beta^{-\ell q_1} \left(\beta^{q_1} -1\right),\\
    \var{G_\ell} = V_\ell &\approx Q_S h_0^{q_2} \beta^{-(\ell-1) q_2}.
\end{align}
\label{eq:asymb_model}
\end{subequations}
Specifically, as a consequence of \eqref{eq:weak_error_model}, the bias model is
\begin{equation}
    B \approx |Q_W| h_0^{q_1} \beta^{-L q_1}.
    \label{eq:bias}
\end{equation}
\end{changes}

Finally, we note that the algorithms presented in this work are iterative.
 We  therefore denote by $\overline M_\ell, \overline G_\ell$ and $\overline V_\ell$ the total number of samples of $G_\ell$ generated in all iterations and their
sample average and sample variance, respectively. Explicitly, we
write\begin{changes}\footnote{\begin{changes}For the variance estimator, one can also use the unbiased
  estimator; by dividing by $\overline M_\ell-1$ instead of $\overline
  M_\ell$. All discussion in this work still applies.\end{changes}}
\begin{subequations}
\begin{align}
    \overline G_\ell &= \frac{1}{\overline M_\ell} \sum_{m=1}^{\overline M_\ell} G_{\ell}(\omega_{\ell,m}),\\
     \label{eq:total_var} \overline V_\ell &= \frac{1}{\overline M_\ell} \sum_{m=1}^{\overline M_\ell} \left( G_{\ell}(\omega_{\ell,m}) - \overline G_\ell \right)^2.
\end{align}
\end{subequations}
\end{changes}

\section{Standard MLMC}\label{s:std} % std.tex
\subsection{Overview}
While minor variations exist among MLMC algorithms listed in \cite{giles08,gs13b,gs13}, we believe that there is sufficient commonality
in them for us to outline here the overarching idea and refer to this collection of methods as the Standard MLMC algorithm or simply SMLMC.
SMLMC solves the problem by iteratively increasing the number of levels of the MLMC hierarchy.
In order to find the optimal number of samples of each level $\ell$, an estimate of the variance $V_\ell$ is needed.
If there were previously generated samples in previous iterations for a level $\ell$, the sample variance $\overline V_\ell$ is used.
Otherwise, an initial fixed number of samples, $\widetilde M$, is generated.
Moreover, in most works, the splitting between bias and statistical error, $\theta$, is chosen to be 0.5.

After running the hierarchy, an estimate of the total error is computed.
To this end, the work \cite{giles08}
approximates the absolute value of the constant, $Q_W$, using a similar expression to the following:
\begin{changes}
\begin{align*}
  |Q_W| \approx \frac{ \max \left(|\overline G_L|, {|\overline
        G_{L-1}|}{\beta^{-q_1}} \right)}{h_0^{q_1} \beta^{-L q_1} (
    \beta^{q_1} -1 )} := \est Q_W .
\end{align*}
\end{changes}
In other words, the absolute value of the constant $Q_W$ is estimated using the samples generated on the last two levels.
Thus, this estimate is only defined for $L\geq2$. Next, the variance of the estimator, $\var{\mathcal A}$, is approximated by
\begin{align*}
  \var{\mathcal{A}} \approx \sum_{\ell=0}^L\frac{\overline V_\ell}{\overline M_\ell} := \est V.
\end{align*}
Finally, a total error estimate can be computed as outlined by \eqref{eq:total-error}
\begin{align}
  \text{Total error estimate} = \est Q_W h_0^{q_1} \beta^{-L q_1} + C_\alpha \sqrt{\est V}.
  \label{eq:std_err_est}
\end{align}
The complete algorithm is outlined in Algorithm \ref{alg:standard}.
\begin{algorithm}[!ht]
\begin{algorithmic}[1]
\Function{StandardMLMC}{$\tol, \est M, \theta$}
    \State Start with $L=0$.
    \Loop
    \State Add new levels to $\{ h_\ell \}_{\ell=0}^L$.
    \State Generate $\est M$ samples for level $L$ and estimate $\overline V_L$.
    \State \begin{minipage}[t]{0.8\textwidth}
            Using sample variance estimates, $\{\overline V_\ell\}_{\ell=0}^L$ from all iterations, and the constant $\theta$, compute optimal number of samples,  $\{M_\ell\}_{\ell=0}^L$,  according to \eqref{eq:optimal_ml}.
            \end{minipage}
    \State Run the hierarchy using the optimal number of samples.
    \State If $L\geq2$ and the total estimate error  \eqref{eq:std_err_est} is less than $\tol$, then \textbf{END}.
    \State Otherwise, set $L = L+1$.
    \EndLoop
\EndFunction
\end{algorithmic}
\caption{}
\label{alg:standard}
\end{algorithm}

Usually all samples from previous iterations are used in the algorithm to run the hierarchy in step 7 to calculate the required quantity of interest.
However, the analysis of the bias and the statistical error of the resulting estimator is difficult and has not been done before, to the best of our knowledge.
\subsection{Accuracy of the parameter estimates}
\label{s:accuracy}
In the standard algorithm, $Q_W$ and the variances $\{V_\ell\}_{\ell=0}^L$ are needed and estimated. In this section, we look at the accuracy
of the estimators for these problem parameters.

We examine the accuracy of the sample variance by computing its squared relative error for $\ell>1$:
\begin{changes}\begin{align*}
\frac{\var{\ \overline V_\ell}}{V_\ell^2} &=
\frac{\left(\overline M_\ell-1\right)^2}{\overline M_\ell^3 V_\ell^2} \left(  \E{\left(G_\ell -\E{G_\ell}\right)^4} -
        \frac{V_\ell^2 (\overline M_\ell - 3)}{\overline M_\ell-1}\right)\\
        &= \frac{\left(\overline M_\ell-1\right)^2}{\overline M_\ell^3} \left( \E{\left(G_\ell -\E{G_\ell}\right)^4} V_\ell^{-2} -
                \frac{\overline M_\ell - 3}{\overline M_\ell-1}\right)\\
        &\approx \frac{\left(\overline M_\ell-1\right)^2}{\overline M_\ell^3} \left( \E{\left(G_\ell -\E{G_\ell}\right)^4} Q_S^{-2} h_\ell^{-2q_2} -
                \frac{\overline M_\ell - 3}{\overline M_\ell-1}\right).
\end{align*}\end{changes}
Unless $\E{\left(G_\ell -\E{G_\ell}\right)^4} \leq C h_\ell^{2q_2},$ for some constant $C>0$,
or $M_\ell$ increases sufficiently fast, the relative error in the estimator $\overline V_\ell$ can become unbounded as $\ell \to \infty$.
Similarly, the relative error of the sample variance at level $\ell=0$ can be shown to be bounded for instance by assuming that the second and fourth central moments of $G_0$ are bounded.

Next, for simplicity, we look at the squared relative error estimate
of $Q_W$ by assuming that it is estimated using samples on a single
level, $L$, only.
\begin{changes}
\begin{align*}
    \frac{\var{\left|\frac{\overline G_L}{h_0^{q_1} \beta^{-L q_1} (  \beta^{q_1}-1)} \right|} }{Q_W^2}
    &= \frac{V_L}{Q_W^2 M_L h_0^{2 q_1} \beta^{-2 L q_1}(\beta^{q_1}-1)^2} \\
    &= \frac{Q_S}{Q_W^2} \cdot \frac{h_0^{q_2} \beta^{-q_2 L}}{Q_W^2 M_L h_0^{2 q_1} \beta^{-2 L q_1}( \beta^{q_1}-1)^2}  \\
    &= \frac{Q_S h_0^{q_2 -2 q_1}}{Q_W^2(\beta^{q_1}-1)^2} \left( \frac{\beta^{L(2q_1-q_2)}}{M_L} \right).
  \end{align*}
  \end{changes}
Observe now that if $q_2<2q_1$ \begin{changes}(as in Example \ref{ex:sde_problem})\end{changes}, then, for the previous relative error estimate to be $\order{1}$, we must have
${M_L \propto \beta^{L(2q_1 -q_2)} \to \infty}$ as $L \to\infty$.
This analysis shows that in some cases, $M_L$ will have to grow to provide an accurate estimate to $Q_W$,
    regardless of the optimal choice of the number of samples outlined in \eqref{eq:optimal_ml}.

\section{Continuation MLMC (CMLMC)}\label{s:adapt} % adapt.tex
In this section we discuss the main contribution of this work, a
continuation MLMC (CMLMC) algorithm that approximates the value
$\E{g(u)}$. We begin in the next subsection by giving an overview of
the general idea of algorithm. Subsequent subsections discuss how to
estimate all the required problem parameters that are necessary for
running the algorithm. CMLMC is listed in
Algorithm~\ref{alg:adaptive}.
\subsection{Overview}
\begin{changes}
The main idea of CMLMC is to solve for $\E{g(u)}$ with a sequence of
decreasing tolerances.
By doing this, CMLMC is able to increasingly improve estimates of several problem
dependent parameters while solving relatively inexpensive problems
corresponding to large tolerances. These parameters estimates are
crucial to optimally distribute computational effort when solving for
the last tolerance, which is the desired tolerance, $\tol$, or smaller.
Moreover, the sequence is built such that the total work of the algorithm is close
to the work of MLMC when solving for the desired tolerance, $\tol$,
assuming all the necessary parameters are known a priori.
To this end, we make the following choice for the sequence of
decreasing tolerances $\tol_i$ for $i=0,1,\ldots$
\end{changes}
\begin{equation*}
    \tol_i = \begin{cases}
        r_1^{i_E-i} r_2^{-1} \tol & i < i_E, \\
        r_2^{i_E-i} r_2^{-1} \tol & i \ge i_E,
        \end{cases}
\end{equation*}
where $r_1 \geq r_2 > 1$. By imposing $\tol_0 =\tol_{\max}$ for some maximum tolerance, we have
\[ i_E = \left\lfloor \frac{-\log(\tol) + \log(r_2) + \log(\tol_{\max})}{\log(r_1)} \right\rfloor,\]

Iterations for which $i \leq i_E$ are meant to obtain increasingly more accurate estimates of the problem parameters.
The iteration $i_E$ solves the problem for the tolerance ${r_2^{-1}\tol}$. Notice that the problem is solved for a slightly smaller tolerance
than the required tolerance $\tol$. This tolerance reduction is to prevent extra unnecessary iterations due to slight variations in estimates of the problem parameters.
This technique improves the overall average running time of the algorithm.
Similarly, iterations $i > i_E$ have tolerances that are even smaller to account for cases in which estimates of the problem parameters are unstable.
The parameters $r_1$ and $r_2$ are chosen \mychanges{such that the total work of the algorithm
is not significantly more than the work of the final hierarchy that solves the problem with the required tolerance, $\tol$}.
For example, if the work of the MLMC estimator is ${\mathcal{O}(\tol^{-2})}$, we choose $r_1=2$ to ensure that
the work of iteration $i$ is roughly four times the work of iteration $i-1$ for iterations for which  ${\tol_i \geq \tol}$.
The choice of $r_2=1.1$, on the other hand, ensures that for iterations for which ${\tol_i < \tol}$, the work of iterations of $i$ is roughly 1.2 times the work of iteration $i-1$.

Consider now the $i$-th iteration of CMLMC and assume that estimates
for ${\mathbf Q := \{q_1, q_2, Q_W, Q_S\}}$ and
$\{V_\ell\}_{\ell=0}^L$ are available from previous iterations; we
will discuss how to obtain these estimate in Section
\ref{ss:param_est}. The $i$-th iteration begins by selecting the
optimal number of levels $L[i]$ that solves the problem for the given
tolerance, $\tol_i$, as follows
\begin{equation}
    L[i] = \text{argmin}_{L_{\min}[i] \leq L \leq L_{\max}[i]} W(\tol_i, L),
    \label{eq:L_optimize}
\end{equation}
where $W(\tol_i,L)$ is defined by \eqref{eq:work_per_L} and depends on all the parameters $\mathbf Q$ and $\{V_\ell\}_{\ell=0}^L$ and $\theta = \theta(L)$ given by
\begin{equation}
\theta = 1 - \frac{| Q_W | h_L^{q_1}}{\tol_i} \mychanges{= 1 - \frac{| Q_W |
  h_0^{q_1} \beta^{-L q_1}}{\tol_i},}
    \label{eq:optimal_theta}
\end{equation}
which comes from enforcing that the bias model \eqref{eq:bias} equals $(1-\theta) \tol_i$.
Moreover, $L_{\min}$ should satisfy ${Q_W h_{L_{\min}}^{q_1} = \tol_i}$ or, since we have $h_\ell = h_0 \,\beta^{-\ell},$
\begin{equation*}
    L_{\min}[i] = \max\left(L[i-1], \frac{q_1 \log(h_0) - \log\left(\frac{\tol_i}{|Q_W|}\right)}{q_1 \log{\beta}} \right),
\end{equation*}
where $L[i-1]$ is the number of levels from the previous iteration.
This ensures that $L$ does not decrease from one iteration to the next, which agrees with our intuition that $L$ increases with $\log\left(\tol_i^{-1}\right)$.
On the other hand, $L_{\max}$ is given by other considerations.
For instance, it could be related to the minimum mesh size imposed by memory or computational restrictions.
More practically, to ensure robustness, $L_{\max}$ can be chosen to be $L_{\min} + L_{\text{inc}}$, for a given fixed integer $L_{\text{inc}}$,
so that $L$ has limited increments from one iteration to the next.
Since only few values of $L$ are considered in the optimization \eqref{eq:L_optimize}, it is easy to find the optimal $L$ by exhaustive search.
The choice \eqref{eq:optimal_theta} implies that the statistical constraint \eqref{eq:var_bound} is relaxed (or tightened) depending on the estimated bias of each hierarchy.
The iteration then continues by \mychanges{running the resulting
  hierarchy}
with the optimal number of samples $\{M_\ell\}_{\ell=0}^L$ according to \eqref{eq:optimal_ml}.
Finally the iteration ends by improving the estimates of the problem parameters $\mathbf Q$ and $\{V_\ell\}_{\ell=0}^L$ as well as the quantity of interest based on the newly available samples as described in Section
\ref{ss:param_est}.

To start CMLMC we compute with an initial, relatively inexpensive,
hierarchy. The purpose of using this initial hierarchy is to obtain
rough estimates of the problem parameters. Such a hierarchy cannot
depend on estimates of problem parameters and should have at least
three levels to allow estimating $\mathbf Q$; \mychanges{these three
  levels are needed to be able to extrapolate (or interpolate) the
  weak error and variance estimates on all MLMC levels}.
The algorithm stops when the total error estimate is below the
required tolerance $\tol$.

\subsection{Parameters estimation}
\label{ss:param_est}
\begin{changes}
In this section, we discuss how to improve estimates of the parameters
$\mathbf Q$ as well as the variances $V_\ell$ based on the generated
samples in all iterations and all levels.
\end{changes}For easier presentation, we will also use the following
notation
\begin{changes}
\begin{align*}
    w_\ell(q_1) &=  h_0^{q_1} \beta^{-\ell q_1} (\beta^{q_1}-1) ,\\
    s_\ell(q_2) &=  h_0^{-q_2} \beta^{\ell q_2}.
  \end{align*}
  \end{changes}
Thus, using the notation above, \eqref{eq:asymb_model} becomes
\begin{subequations}
\begin{align}
    \E{G_\ell} &\approx Q_W w_\ell(q_1),\\
    \var{G_\ell} = V_\ell &\approx Q_S s_\ell^{-1}(q_2).
    \label{eq:vl_asymb_model_2}
\end{align}
\label{eq:asymb_model_2}
\end{subequations}

\subsubsection{Estimating variances $V_\ell$}
 \label{ss:var_est}
 We first assume that we have estimates of $q_1$, $q_2$, $Q_W$ and $Q_S$ and discuss estimating the variances, $\{V_\ell\}_{\ell=0}^L$,
 and the total statistical error after computing with a given hierarchy. Estimating $q_1, q_2, Q_W$ and $Q_S$ is discussed in the next subsection.

Usually the variances $\{V_\ell\}_{\ell=0}^L$ are estimated by using the sample variance estimator \eqref{eq:total_var}
to estimate the statistical error as well as the optimal number of samples $\{M_\ell\}_{\ell=0}^L$.
However, sometimes there are too few samples in a given level to give a corresponding accurate variance estimate.
This is specially acute  on the deepest levels,
and unlike the standard MLMC algorithm, we do not impose a minimum number of samples across levels
to obtain a stable estimate of the sample variance.
Recalling that we have the variance model \eqref{eq:vl_asymb_model_2} at our disposal,
we can use this model to estimate the variance at all levels $\ell > 0$.
However, the model \eqref{eq:vl_asymb_model_2} is only accurate asymptotically.
We can use the generated samples on each level to locally improve the accuracy of the $V_\ell$ estimates.
To this end, we use a Bayesian setting \cite{sivia1996data}.

We assume that $G_\ell$ follows a normal distribution with mean
$\mu_\ell$ and precision $\lambda_\ell$ (precision is simply the inverse of the variance).
To simplify the computation, we choose a normal-gamma prior on $(\mu_\ell, \lambda_\ell)$ -- the conjugate prior of the normal likelihood.
The resulting posterior  probability density function (pdf) is also a
normal-gamma distribution function.
\begin{changes}
We choose the parameters $(\widehat \mu_\ell,\kappa_0, 0.5 + \widehat \lambda_{\ell}\kappa_1, \kappa_1)$ for the normal-gamma prior, such that it is maximized at $\widehat \mu_{\ell}$ and $\widehat \lambda_{\ell}$.\end{changes}
The parameter $\widehat \mu_\ell$ and $\widehat \lambda_\ell$ serve as initial guesses for $\mu_\ell$ and $\lambda_\ell$, respectively.
Moreover, $\kappa_0$ and $\kappa_1$ are positive constants that model our certainty in those respective guesses.
We use the assumed models of the weak and strong errors \eqref{eq:asymb_model_2} to give the initial guesses
\begin{subequations}
\begin{align}
    \widehat \mu_{\ell} &= Q_W w_\ell(q_1),\\
    \label{eq:vl_initial_guess}
    \widehat \lambda_{\ell} &=  Q_S^{-1} s_\ell(q_2).
\end{align}
\end{subequations}
As mentioned, the posterior pdf is also a normal-gamma with parameters
$(\Upsilon_{1,\ell},\Upsilon_{2,\ell},\Upsilon_{3,\ell},\Upsilon_{4,\ell})$ and it is maximized at
\begin{changes}$\left(\Upsilon_{1,\ell}, \frac{\Upsilon_{3,\ell}-0.5}{\Upsilon_{4,\ell}}\right)$\end{changes}. Specifically
\begin{align*}
    \Upsilon_{3,\ell}  &= 0.5 + \kappa_1 \widehat \lambda_\ell + \frac{\overline M_\ell}{2},\\
    \Upsilon_{4,\ell}  &= \kappa_1 + \frac{1}{2} \left( \sum_{m=1}^{\overline M_\ell} \left( G_{\ell, m} -\overline G_\ell \right)^2 \right)
                + \frac{\kappa_0 \overline M_\ell (\overline G_\ell -\widehat \mu_\ell)^2}{2 (\kappa_0 + \overline M_\ell)}.
\end{align*}
As such, we use the following estimate of the variance $V_\ell$ for $\ell>0$
\begin{equation}
    V_\ell \approx \frac{\Upsilon_{4,\ell}}{\Upsilon_{3,\ell}-0.5}.
    \label{eq:vl_bayes}
\end{equation}
Estimating the variance at the coarsest mesh, $V_0$, can be done using the sample variance.
The number of samples on the coarsest level, $M_0$, is usually large enough to produce a stable and accurate estimate.
Using these estimates and the bias estimate \eqref{eq:bias}, the total error can be estimated as \eqref{eq:total-error}.

\subsubsection{Estimating $\mathbf Q$}
To incorporate prior knowledge on $q_1$ and $q_2$ including initial guesses and the relation $q_2 \leq 2q_1$,
we again follow a Bayesian setting to estimate these parameters and assume that $G_\ell$ follows a Gaussian distribution with mean
${Q_W w_\ell(q_1)}$ and variance ${Q_S s_\ell^{-1}(q_2)}$.
In what follows, $\ell_0$ is a non-negative integer. With these assumptions, the corresponding likelihood is
\begin{align}
    \mathcal{L} = \left(\prod_{\ell=\ell_0}^L \left(2 \pi Q_S s_\ell^{-1}(q_2)\right)^\frac{-\overline M_\ell}{2} \right)
    \exp\left( -\frac{1}{2 Q_S} \sum_{\ell=\ell_0}^L s_\ell(q_2) \sum_{m=1}^{\overline M_\ell} \left( G_{\ell,m} - Q_W w_\ell(q_1) \right)^2 \right).
    \label{eq:likelihood}
\end{align}
Assuming a improper prior on $Q_W$ and $Q_S$ and maximizing the resulting posterior pdf with respect to $Q_W$ and $Q_S$ gives the following weighted least-squares solution:
\begin{subequations}
\begin{align}
   \label{eq:qw_mean} Q_W^* &= \left( \sum_{\ell=\ell_0}^L \overline M_\ell w_\ell^2(q_1) s_\ell(q_2)  \right)^{-1} \sum_{\ell=\ell_0}^L
        w_\ell(q_1) s_\ell(q_2) \overline M_\ell \overline G_\ell ,\\
    Q_S^* &= \left( \sum_{\ell=\ell_0}^L \overline M_\ell \right)^{-1} \sum_{\ell=\ell_0}^L s_\ell(q_2)
                \sum_{m=1}^{\overline M_\ell} \left(G_{\ell,m} - Q_W w_\ell(q_1) \right)^2.
\end{align}
\label{eq:lsq-est}
\end{subequations}
We can substitute the previous expressions for $Q_W$ and $Q_S$ in \eqref{eq:likelihood} to obtain a likelihood in terms of $q_1$ and $q_2$.
Denoting $\overline M = \sum_{\ell=\ell_0}^L \overline M_\ell$, we write
\begin{align*}
    \mathcal{L}(q_1, q_2) = \exp\left(-\frac{\overline M}{2}\right)
        \left( \sum_{\ell=\ell_0}^L \sum_{m=0}^{\overline M_\ell} s_\ell(q_2) G^2_{\ell,m} -
                \frac{\left(\sum_{\ell=\ell_0}^L s_\ell(q_2) w_\ell(q_1) \overline M_\ell \overline G_{\ell}\right)^2}{\sum_{\ell=\ell_0}^L \overline M_\ell w_\ell(q_1)^2
                s_\ell(q_2)
                    } \right)^{-\frac{\overline M}{2}} .
\end{align*}
We can then assume a prior on $q_1$ and $q_2$.
However, remember that $q_2 \leq 2 q_1$, and $q_1 > 0$.  As such, we introduce the unconstrained parameters
${x_0(q_1) = \log(q_1) \in \rset}$ and ${x_1(q_1,q_2) = \log(2 q_1 - q_2) \in \rset}$ and assume a Gaussian prior on them
\begin{align*} \rho_{\text{prior}}(q_1,q_2) = \frac{1}{2 \pi \sqrt{\sigma_0^2 \sigma_1^2}} \exp \left(-\frac{(x_0(q_1) - \widehat x_0)^2}{2 \sigma_0^2} -
        \frac{(x_1(q_1,q_2) - \widehat x_1)^2}{2 \sigma_1^2}\right).
\end{align*}
Here, $\widehat x_0$ and $\widehat x_1$ represent our initial guesses of $x_0$ and $x_1$, respectively, which we can obtain from a rough analysis of the problem.
Moreover, $\sigma_1$ and $\sigma_2$ model our confidence in those guesses. The more accurate our initial guesses are, the faster the algorithm converges.
Finally, we numerically maximize the log of the posterior pdf with respect to $(x_0,x_1) \in \rset^2$ using a suitable numerical optimization algorithm.
For robustness, we choose $\ell_0=1$ to estimate $q_1$ and $q_2$. In other words we include samples from all levels $\ell>0$ for this estimation.

 \mychanges{Given estimates of $q_1$ and $q_2$, we can use the least-squares
 estimates $Q_W^*$ and $Q_S^*$ in \eqref{eq:lsq-est} as estimates of
 $Q_W$ and $Q_S$, respectively}.
 However, usually not all levels follow the assumed asymptotic models \eqref{eq:asymb_model} and as such special care must be taken to choose
 $\ell_0$ in these estimates.
 The parameter $Q_W$ must be accurate on deeper levels since it is used to compute the bias \eqref{eq:bias}. Similarly, $Q_S$ must be accurate on deeper
 levels where not many samples are available and the variance estimate \eqref{eq:vl_bayes} is mainly determined by the initial guess
 \eqref{eq:vl_initial_guess}.
 For these reasons, \mychanges{when computing $Q_W^*$ and $Q_S^*$}, we choose $\ell_0=\max(1,L-\mathfrak{L})$ in \eqref{eq:lsq-est}
 for some positive integer $\mathfrak{L}$ that denotes the maximum number of levels use to compute the estimates.
\mychanges{Finally, since $Q_W$ has an improper prior, its posterior is also a
 Gaussian with mean $Q_W^*$ and variance
 \[V_W := \sum_{\ell=\ell_0}^L \frac{Q_S}{\overline M w_\ell^2(q_2) s_\ell(q_1)}. \]
 Motivated by the accuracy analysis of the $Q_W$ estimate in Section
 \ref{s:accuracy},
 we use a worst-case estimate of $Q_W$
 instead of simply using the estimate $Q_W^*$ in \eqref{eq:qw_mean},
 The worst-case estimate is produced by adding the maximum sampling error with $1-\alpha$ confidence,
 namely $C_\alpha \sqrt{V_W}$,
 multiplied by the sign of $Q_W^*$. In
 other words, our
 estimate of $Q_W$ is $Q_W^* + \text{sign}(Q_W^*) C_\alpha \sqrt{V_W}$.}

\begin{changes}
\subsection{Algorithm parameters}
Table \ref{tbl:adapt_parameters} summarizes the parameters that control the
CMLMC algorithm. Some of these parameters need to be suitably chosen
for the specific problem.
However, while there might be optimal values for these
parameters to minimize the average running time, it is our experience
that reasonable values of these parameters are enough to get average
running times that are near-optimal.
In fact, similar results to those that we show
Section~\ref{s:res} were obtained with variations of $\kappa_1$ and $\kappa_2$; namely $\kappa_1 = \kappa_2 \in \{0.05, 0.1, 0.2\}$.
\end{changes}

\begin{algorithm}
\begin{algorithmic}[1]
\Function{CMLMC}{Parameters summarized in Table \ref{tbl:adapt_parameters}}
    \State Compute with an initial hierarchy.
    \State \begin{minipage}[t]{0.8\textwidth}
            Estimate problem parameters $\left\{ V_\ell\right\}_{\ell=0}^L, Q_S, Q_W, q_1$ and $q_2$
               according to section \ref{ss:param_est}.
            \end{minipage}
    \State Set $i = 0$.
    \Repeat
        \State Find $L$ according to \eqref{eq:L_optimize}.
       \State \mychanges{Add new levels to $\{ h_\ell \}_{\ell=0}^L$}.
        \State \begin{minipage}[t]{0.8\textwidth}
                    Using the variance estimates \eqref{eq:vl_bayes} and $\theta$ from \eqref{eq:optimal_theta}, compute the optimal number of
                    samples according to \eqref{eq:optimal_ml}.
                \end{minipage}
       \State \mychanges{Run} the resulting hierarchy using the optimal number of samples.
       \State \begin{minipage}[t]{0.8\textwidth}
                Estimate problem parameters, $\left\{ V_\ell\right\}_{\ell=0}^L, Q_S, Q_W, q_1$ and $q_2$,
                according to section \ref{ss:param_est}.
                \end{minipage}
       \State Estimate the total error according to \eqref{eq:total-error}.
       \State Set $i = i+1$
        \Until{$i > i_E$ and the total error estimate is less than $\tol$}
\EndFunction
\end{algorithmic}
\caption{}
\label{alg:adaptive}
\end{algorithm}

\begin{table}
\centering
\begin{tabular}{r|p{3.5in}}
\hline
Parameter & Purpose \\
\hline
$\widehat x_0, \widehat x_1, \sigma_0$ and $\sigma_1$ & Parameters to model the initial guess of $q_1$ and $q_2$ and the confidence in those estimates.\\
$\kappa_0$ and $\kappa_1$ &  The confidence in the weak and strong error models, respectively. \\
$\tol_{\max}$ & The maximum tolerance with which to start the algorithm. \\
$r_1$ and $r_2$  & Controls the computational burden to calibrate the problem parameters compared to the one taken to solve the problem. \\
Initial hierarchy & The initial hierarchy to start the algorithm. Must be relatively inexpensive and has at least three levels.  \\
$L_{\text{inc}}$ & Maximum number of values to consider when optimizing for $L$. \\
$\mathfrak{L}$ & Maximum number of levels used to compute estimates of $Q_W$ and $Q_S$.  \\
$C_\alpha$ & Parameter related to the confidence in the statistical constraint.\\
\hline
\end{tabular}
\caption{Summary of parameters in CMLMC}
\label{tbl:adapt_parameters}
\end{table}

\section{Numerical Tests}\label{s:res} % res.tex
In this section, we first introduce the test problems.
 We then describe several implementation details and finish by presenting the actual
numerical results.

\subsection{Test Problems}  \label{ss:examples}
We look at three test problems: the first two are \mychanges{based on} PDEs with random inputs and
the last one is \mychanges{based on} an It\^o SDE.

\subsubsection{Ex.1}
 This problem is based on Example~\ref{ex:spde_problem} in Section~\ref{sec:hier_intro} with some particular choices that satisfy the assumptions therein.
 First, we choose ${\mathcal{D} = [0,1]^3}$ and assume that the forcing is
 \[ f(\mathbf x;\omega) = f_0 + \widehat f \sum_{i=0}^K \sum_{j=0}^K \sum_{k=0}^K \Phi_{ijk}(\mathbf x) Z_{ijk},\]
 where
 \[ \Phi_{ijk}(\mathbf x) = \sqrt{\lambda_i \lambda_j \lambda_k} \phi_i(x_1)\phi_j(x_2)\phi_k(x_3), \]
and
\begin{align*}
    \phi_i(x) &=
\begin{cases}
    \cos\left( \frac{5 \Lambda i}{2} \pi x \right) & i \text{ is even}, \\
    \sin\left( \frac{5 \Lambda (i+1)}{2} \pi x \right) & i \text{ is odd},
\end{cases}, \\
\lambda_i &=  \left( 2 \pi \right)^\frac{7}{6}\Lambda^\frac{11}{6}
\begin{cases}
    \frac{1}{2}  & i=0, \\
    \exp\left( - 2\left(  \frac{\Lambda i}{4} \pi \right)^2 \right) & i \text{ is even}, \\
    \exp\left( - 2\left(  \frac{\Lambda (i+1)}{4} \pi \right)^2 \right) & i \text{ is odd},
\end{cases}
\end{align*}
for given $\Lambda > 0$, and positive integer $K$ and $\mathbf Z = \lbrace Z_{ijk} \rbrace$ a set of $(K+1)^3$ i.i.d. standard normal random variables.
Moreover, we choose the diffusion coefficient to be a function of two random variables as follows:
\begin{align*}
    a(\mathbf x; \omega) &= a_0 + \exp \Big(4 Y_1 \Phi_{121}(\mathbf x) + 40 Y_2 \Phi_{877}(\mathbf x)\Big).
\end{align*}
Here, $\mathbf Y = \lbrace Y_1, Y_2 \rbrace$ is a set of i.i.d. normal Gaussian random variables, also independent of $\mathbf Z$.
Finally we make the following choice for the quantity of interest,
$g$:
\begin{changes}\[g = \left( 2 \pi \sigma^2 \right)^\frac{-3}{2} \int_\mathcal{D} \exp
\left( - \frac{ \| \mathbf x - \mathbf x_0 \|^2_2}{2 \sigma^2} \right)
u(\mathbf x) d\mathbf x, \]\end{changes}
and select the parameters $a_0=0.01, f_0=50, \widehat f=10, \Lambda = \frac{0.2}{\sqrt{2}}, K=10, \sigma^2=0.02622863$ and
${\mathbf x_0 = \left[0.5026695,0.26042876,0.62141498\right]}$.
Since the diffusion coefficient, $a$, is independent of the forcing, $f$, a reference solution can be calculated to sufficient accuracy
by scaling and taking expectation of the weak form with respect to $\mathbf Z$ to
obtain a formula with constant forcing for the conditional expectation with respect to $Y$.
\mychanges{We then use stochastic collocation \cite{bnt2010} with 11
  Hermite quadrature points in each direction (thus totaling 121 points)
and a Finite Difference method with
  centered differences and 128 equally spaced points in each dimension to produce the reference value $\E{g}$.
Using this method, the reference value $1.6026$  was computed with an error estimate of $10^{-4}$.}

\subsubsection{Ex.2}
The second example is a slight variation of the first.
First, we choose the following diffusion coefficient instead:
\begin{align*}
    a(\mathbf x; \omega) &= a_0 + \exp \Big(Y_1 \phi_{121}(\mathbf x) + Y_2 \phi_{877}(\mathbf x)\Big).
\end{align*}
Moreover, in this example $\mathbf Y$ is a set of two i.i.d. uniform random variables in the range $[-1,1]$, again independent of $\mathbf Z$.
We also make the following choice for the quantity of interest $g$
\begin{changes}\[g = 100 \left( 2 \pi \sigma^2 \right)^\frac{-3}{2} \int_\mathcal{D} \exp \left( - \frac{ \| \mathbf x -
        \mathbf x_0 \|^2_2}{2 \sigma^2} \right) u(\mathbf x) d\mathbf x, \]\end{changes}
and select the parameters $a_0=1, f_0=1, \widehat f=1,\Lambda = 0.2,
K=10, \sigma^2=0.01194691$ and ${\mathbf x_0 =
  \left[0.62482261,0.45530923,0.49862328\right]}$.
\mychanges{We use the same method as in \textbf{Ex.1} to compute the reference
solution, except that in this case, we use Legendre quadrature in the
stochastic collocation method, instead of Hermite quadrature.}
The computed reference solution $\E{g}$ in this case is $2.3627$ with an error estimate of $10^{-4}$.

\subsubsection{Ex.3}
    The third example is a one-dimensional geometric Brownian motion based on Example~\ref{ex:sde_problem}. We make the following choices:
    \begin{align*}
        T &=1, \\
        a(t,u) &= 0.05 u,\\
        b(t,u) &= 0.2 u, \\
        g(u) &= 10 \max(u(1)-1,0).
    \end{align*}
    The exact solution can be computed using a standard change of variables and It\^o's formula.
    For the selected parameters, the solution is $\E{g} = 1.04505835721856$.

\subsection{Implementation and Runs}
All the algorithms mentioned in this work were implemented using the C programming language,
with the goal that the software be as optimal as possible, while maintaining generality.

For implementing the solver for the PDE test problems ({\bf Ex.1} and {\bf Ex.2}), we use
PetIGA~\cite{petiga,Collier2013}. While the primary intent of this framework is to
provide high-performance B-spline-based finite element
discretizations, it is also useful for applications where the domain
is topologically square and subject to uniform refinements. As its
name suggests, PetIGA is designed to tightly couple to
PETSc~\cite{petsc-web-page}. The
framework can be thought of as an extension of the PETSc library, which
provides methods for assembling matrices and vectors related to the discretization of
integral equations.

In our PDE numerical tests ({\bf Ex.1} and {\bf Ex.2}), we use a standard trilinear basis to discretize
the weak form of the model problem, integrating with eight
quadrature points. We also generate results for two linear solvers that PETSc provides an interface to.
The first solver is an iterative GMRES solver that solves a linear system in almost linear time with respect to the number of degrees of freedom
    for the mesh sizes of interest; in other words $\tilde \gamma=1$ in this case.
The second solver we tried is a direct one, called MUMPS~\cite{Amestoy2001,Amestoy2006}. For the mesh sizes of interest,
    the running time of MUMPS varies from quadratic to linear in the total number of degrees of freedom. The best fit turns out to be $\tilde \gamma=1.5$ in the case.

From \mychanges{\cite[Theorem 2.3]{tsgu13}}, the complexity rate for all the examples is expected to be
${\Order{TOL^{-s_1} \log(\tol)^{s_2}}}$, where ${s_1}$ and ${s_2}$ depend on $q_1, q_2$ and $d\gamma$.
These and other problem parameters are summarized in Table~\ref{tbl:problem_params} for the different examples.
\begin{table}
\centering
\begin{tabular}{r|cccccc}
\hline
& $d$ & $\tilde \gamma$ & $q_1$ & $q_2$ & $s_1$ & $s_2$ \\
\hline
{\bf Ex.1} and {\bf Ex.2} with GMRES solver & 3 & 1 & 2 & 4 & 2 & 0 \\
{\bf Ex.1} and {\bf Ex.2} with MUMPS solver & 3 & 1.5 & 2 & 4 & 2.25  & 0\\
{\bf Ex.3} & 1 & 1 & 1 & 1 & 2  & 2 \\
\hline
\end{tabular}
\caption{Summary of problem parameters}
\label{tbl:problem_params}
\end{table}
We run each algorithm $100$ times \begin{changes}for each tolerance\end{changes} and show in plots in the next section the medians with vertical bars spanning from the $5\%$ percentile
to the $95\%$ percentile.
Finally, all results were generated on the same machine with $52$ gigabytes of memory to ensure that no overhead is introduced due to hard disk access
during swapping that could occur when solving the three-dimensional PDEs with a fine mesh.

In order to compare CMLMC to SMLMC, and since the latter does not include a step to fit $q_1$ and $q_2$,
we assume that these parameters are both known as discussed in  Example~\ref{ex:spde_problem} and Example~\ref{ex:sde_problem}.
Moreover, we use the parameters listed in Table~\ref{tbl:alg_params}.

\subsection{Results}
\begin{figure}[ht]
\begin{center}
\includegraphics[scale=0.45]{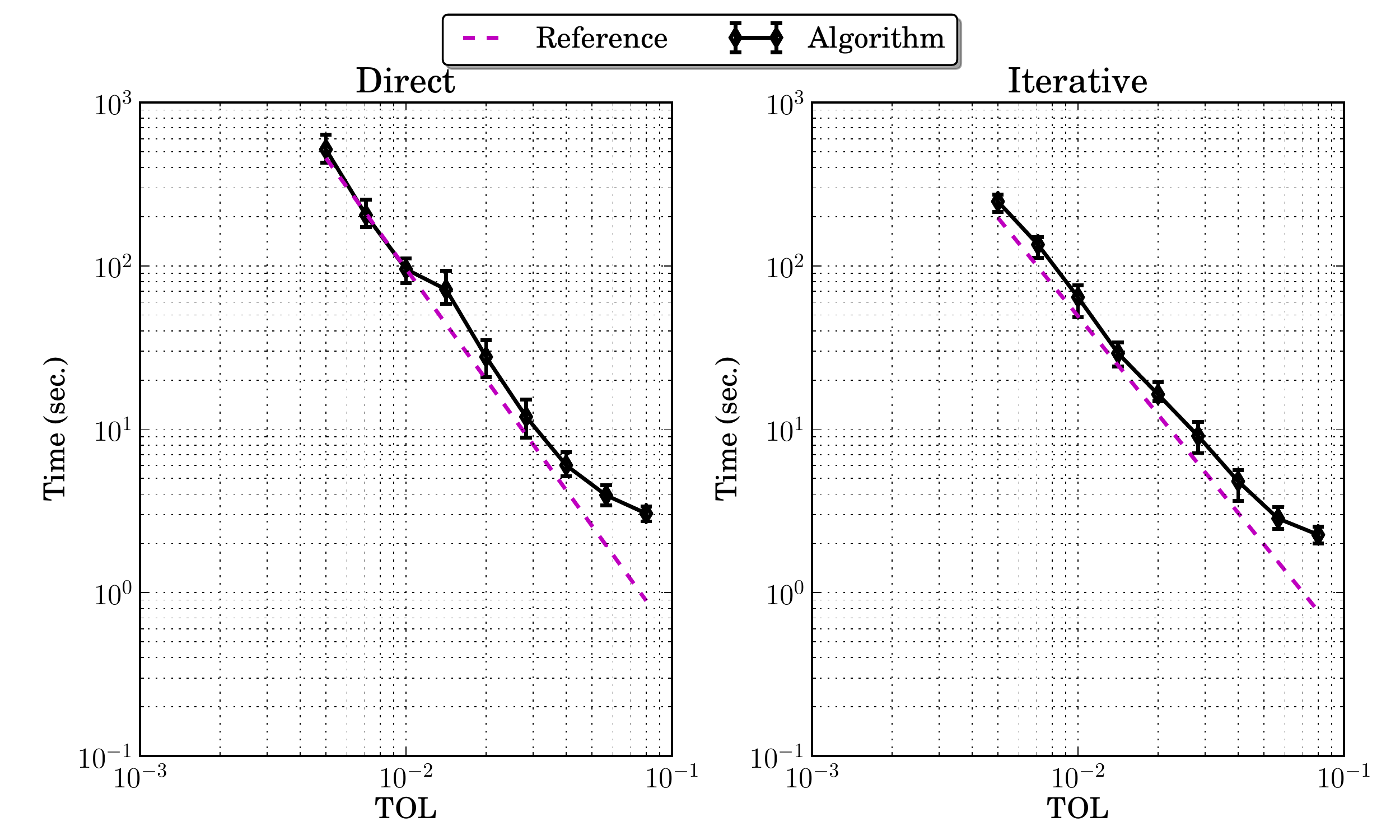}
\includegraphics[scale=0.45]{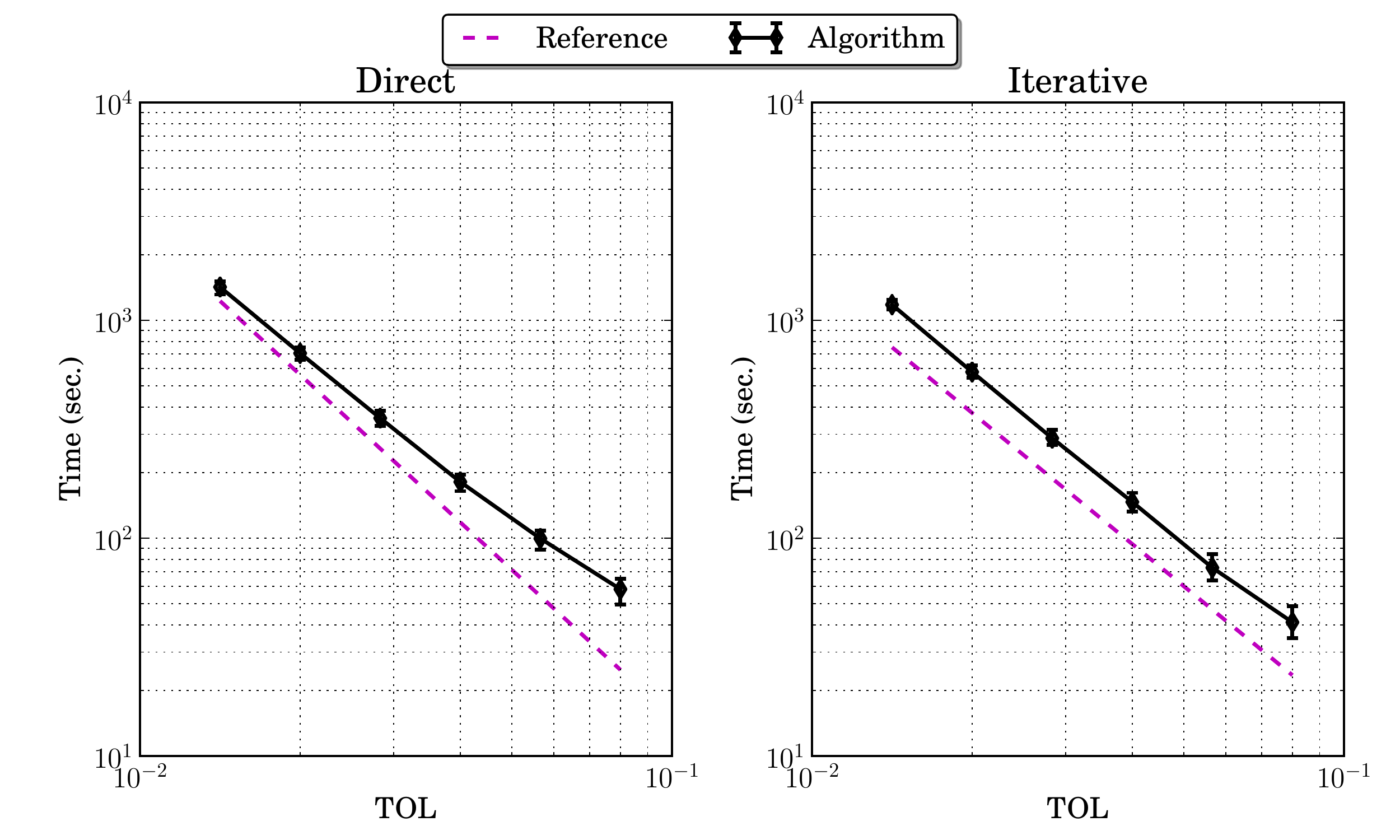}
\includegraphics[scale=0.35]{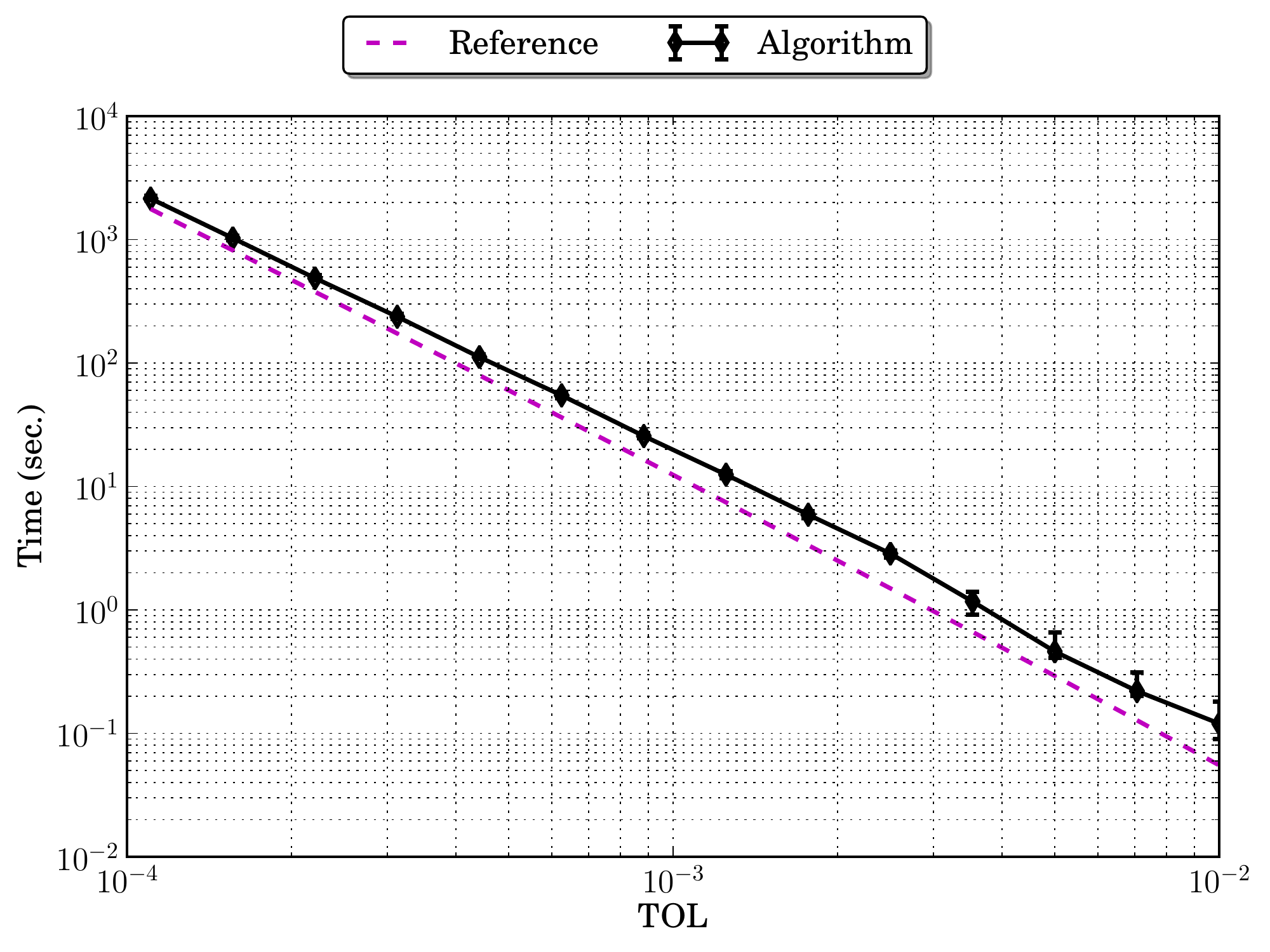}
\end{center}
\caption{\mychanges{From top: {\bf Ex.1, Ex.2, Ex.3}. These plots show
    the total running time of CMLMC. The reference dashed lines are
    $\Order{\tol^{-s_1}\log(\tol)^{s_2}}$ as summarized in
    Table~\ref{tbl:problem_params}. Notice that, asymptotically, the
    total running times seem to follow the expected rates. This shows
    that the algorithm, in our examples, has the same complexity as
    the theoretical work (c.f. \cite[Theorem 2.3]{tsgu13}) of the last
    iteration where we effectively solve the problem with the required
    tolerance requirements.}}
\label{fig:runtime-rate}
\end{figure}

\begin{figure}[ht]
\begin{center}
\includegraphics[scale=0.45]{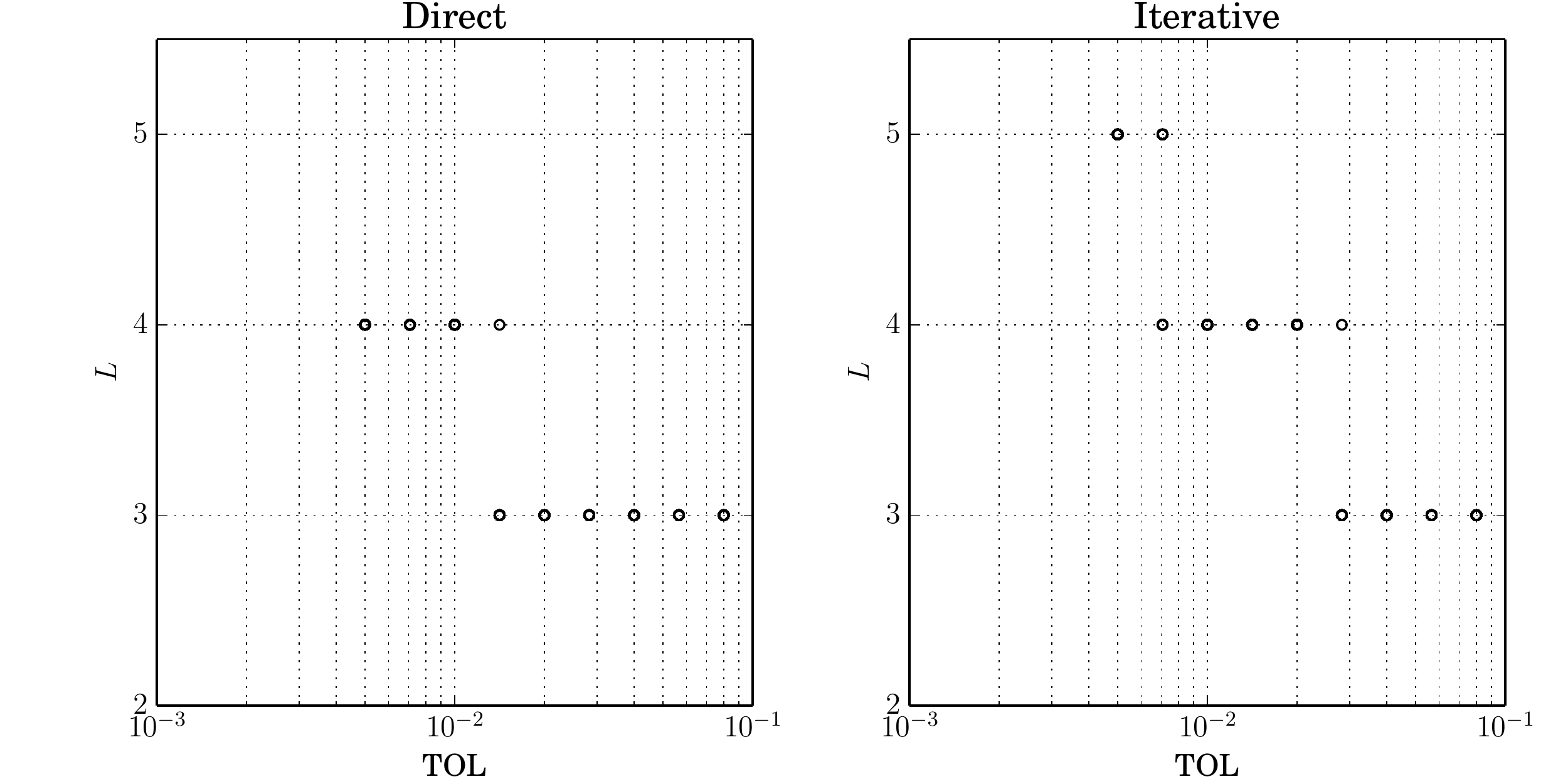}
\includegraphics[scale=0.45]{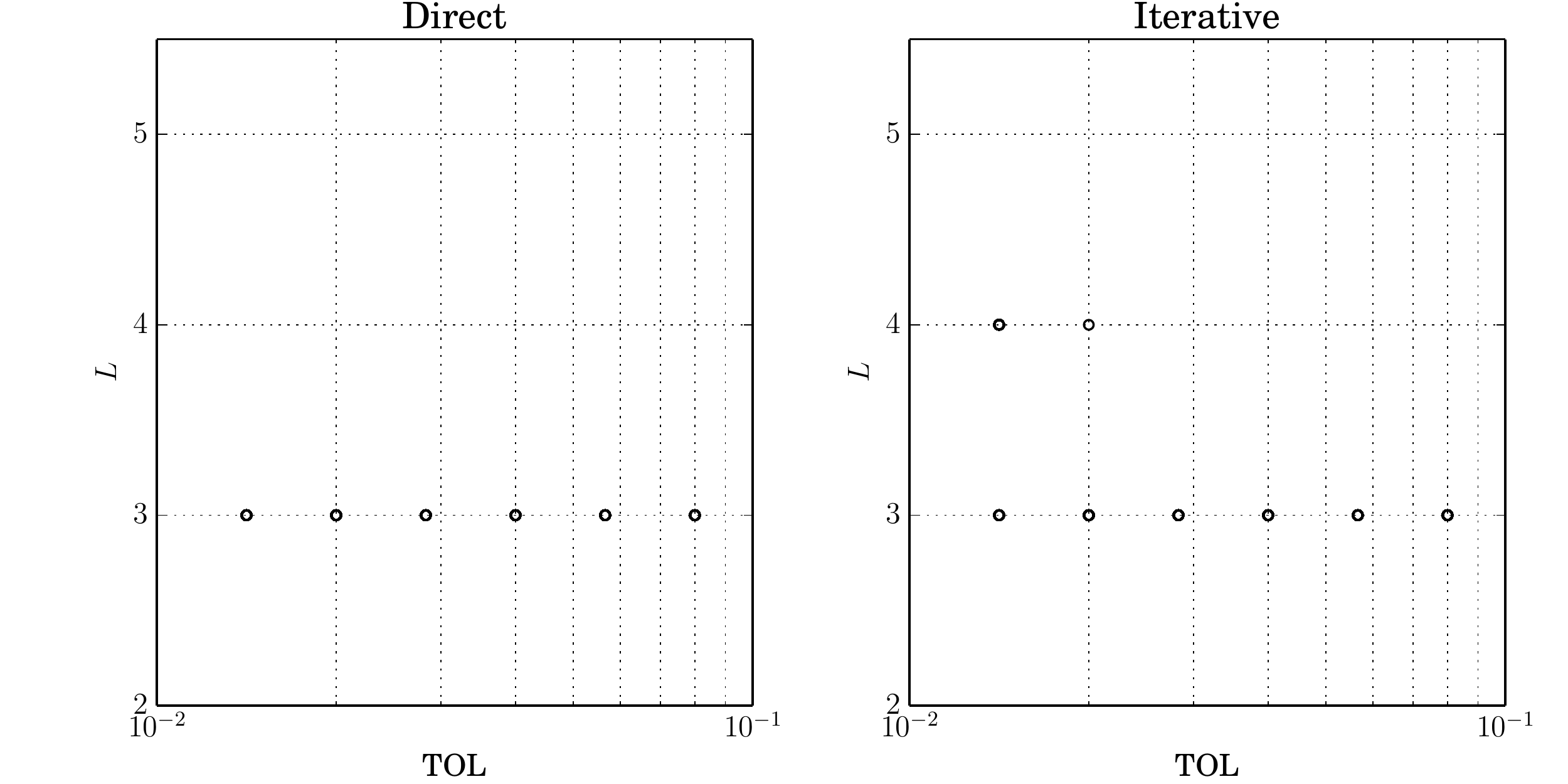}
\includegraphics[scale=0.35]{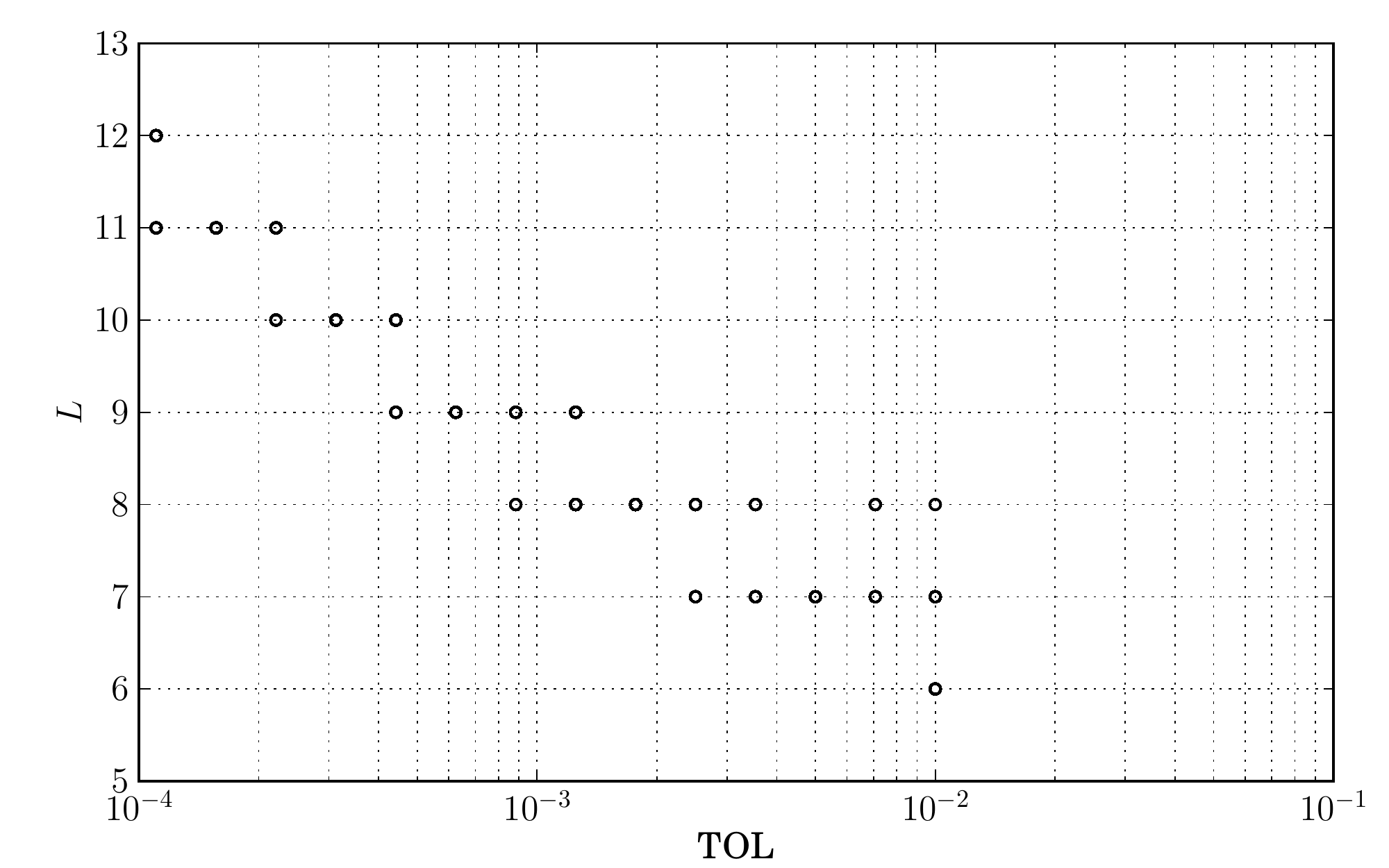}
\end{center}
\caption{\mychanges{From top: {\bf Ex.1, Ex.2, Ex.3}. These plots show
    the number of levels, $L$, for different tolerances, as produced
    in the last iteration of CMLMC. Here, it is clear that $L$ depends
    on the particular realization. However, the relation between $L$
    and $\log(\tol^{-1})$ looks reasonably linear, as expected.
  Note that in {\bf Ex.2}, $L$ does no exhibit significant variations.
This is because, for the tolerances considered, $L=3$ already
satisfies the bias constraint.} }
\label{fig:L-vs-tol}
\end{figure}

\begin{table}
\centering
\begin{tabular}{r|p{1.75in}|p{1.75in}}
\hline
Parameter & Value for PDE examples ({\bf Ex.1} and {\bf Ex.2}) & Value for SDE example ({\bf Ex.3})\\
\hline
$h_0$    & 1/4 for {\bf Ex.1}, 1/8 for {\bf Ex.2} & 1 \\
$\beta$  & 2 & 2 \\
$\kappa_0$ and $\kappa_1$ &  $0.1$ for both & $0.1$ for both\\
$\tol_{\max}$ & 0.5 & 0.1\\
$r_1$ and $r_2$  & $2$ and $1.1$, respectively & $2$ and $1.1$, respectively \\
Initial hierarchy & $L=2$ and $h_\ell = \{4,6,8\}$ and $M_\ell = 10$ for all $\ell$. & $L=2$ and $h_\ell = \{1,2,4\}$ and $M_\ell = 10$ for all $\ell$.  \\
$L_{\text{inc}}$ & 2 & 2  \\
$\mathfrak{L}$ & 3 & 5  \\
$C_\alpha$ & 2 & 2  \\
\hline
\end{tabular}
\caption{Summary of parameters values to used in numerical tests}
\label{tbl:alg_params}
\end{table}
Figure~\ref{fig:runtime-rate} shows that the running time of CMLMC
follows the expected complexity rates ${\Order{\tol^{s_1}
    \log(\tol)^{s_2}}}$ as summarized in
Table~\ref{tbl:problem_params}. \mychanges{Notice that the running
  time in this and all figures that we present in this work include
  the time necessary to sample the underlying stochastic solution and
  the time to do the necessary computation to estimate the problem
  parameters. However, the computational complexity of calculating of
  problem parameters is largely dominated by the computational
  complexity of sampling the approximate solution to the differential
  equations at hand. Indeed, the computations described in Section
  \ref{ss:param_est} are inexpensive post-processing calculations of
  the these samples. Moreover, our results show that the algorithm has
  the same complexity as the theoretical work (c.f. \cite[Theorem
  2.3]{tsgu13}) of the last iteration where we effectively solve the
  problem with the required tolerance requirements.}
Next, Figure~\ref{fig:L-vs-tol} shows the number of levels, $L$, in
the last iteration of CMLMC for different tolerances. As expected,
even though $L$ depends on the particular realization, it is well
approximated by a linear function of $\log(\tol^{-1}).$

\begin{figure}[ht]
\begin{center}
\includegraphics[scale=0.45]{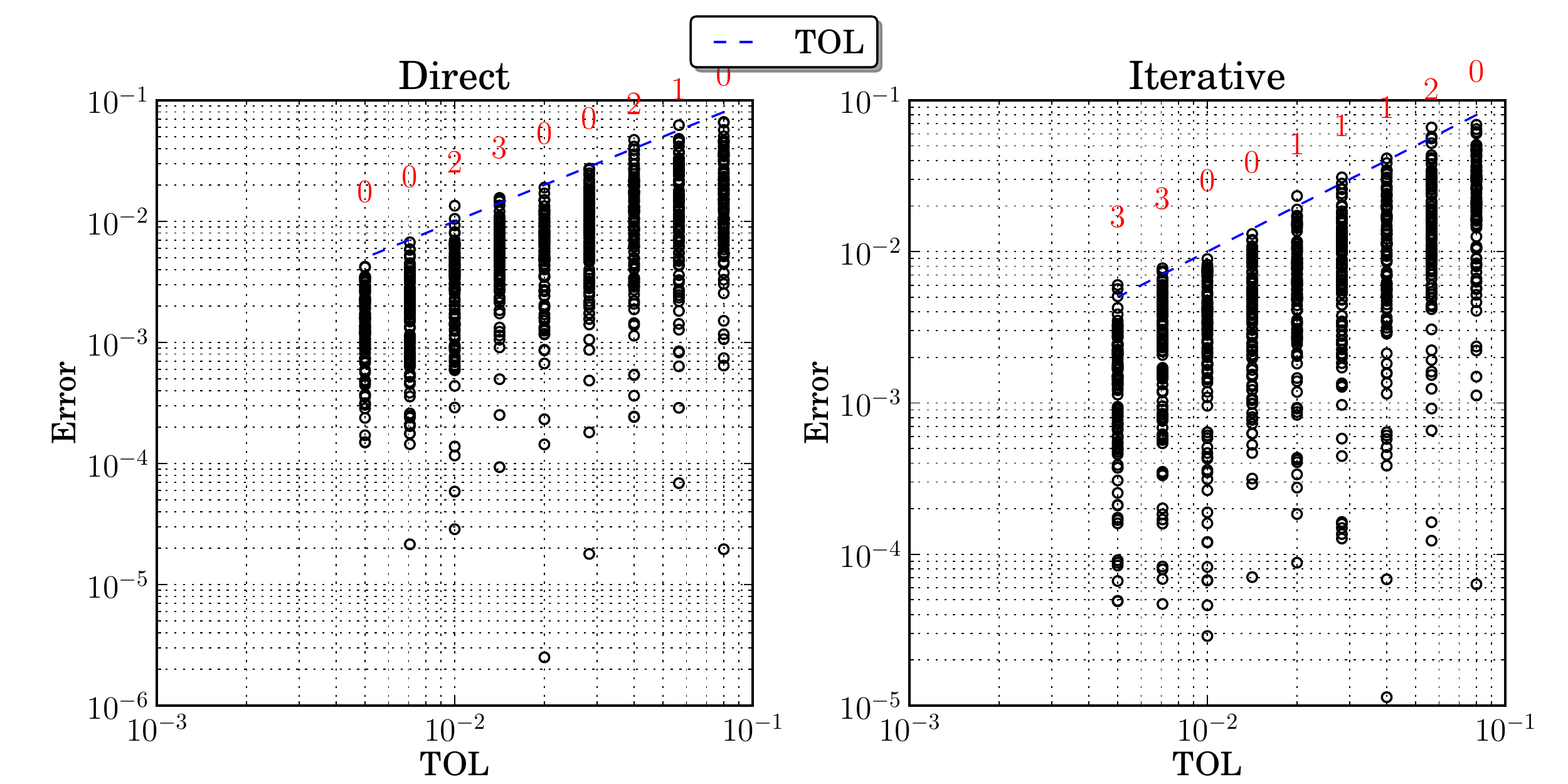}
\includegraphics[scale=0.45]{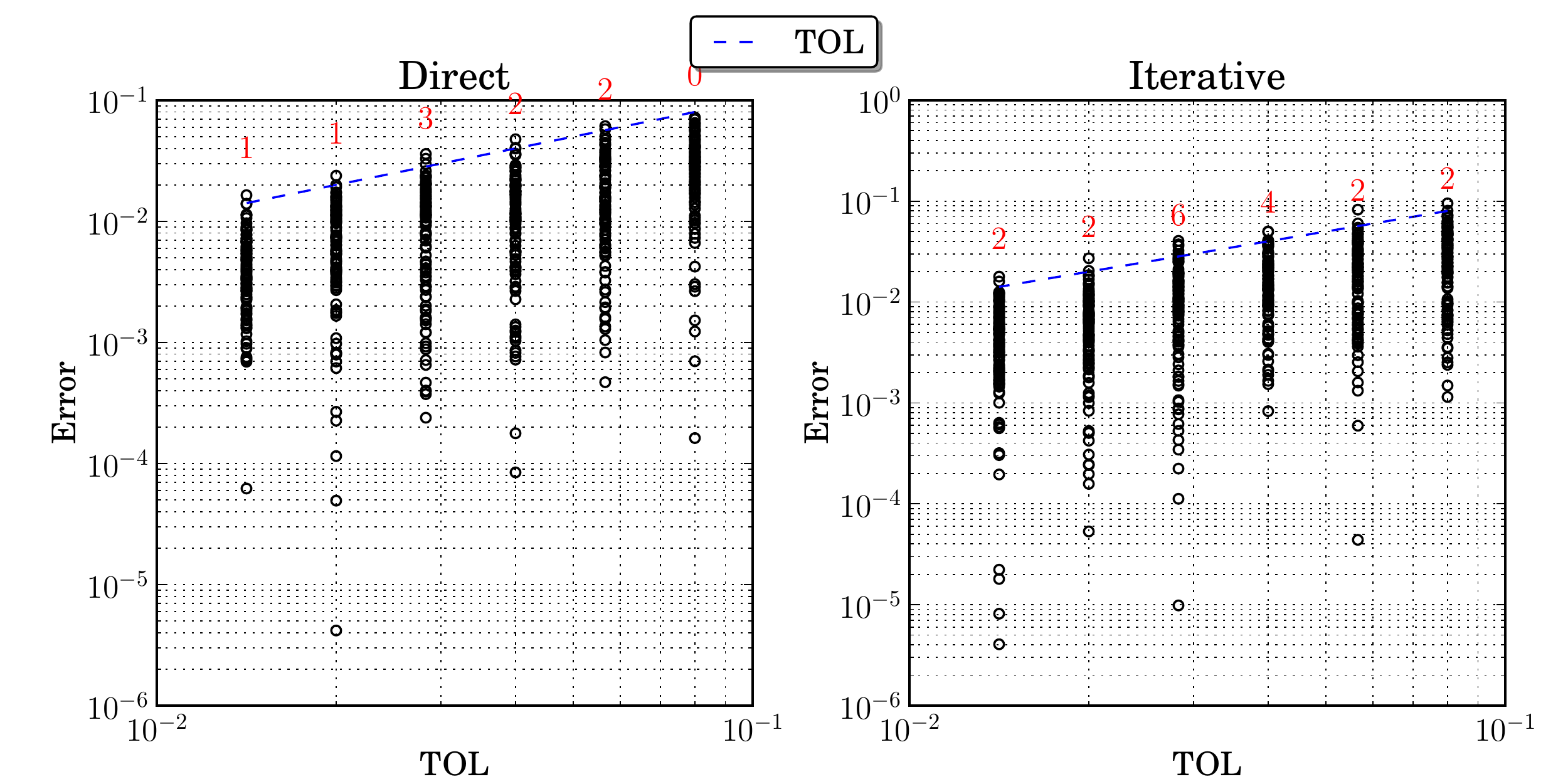}
\includegraphics[scale=0.35]{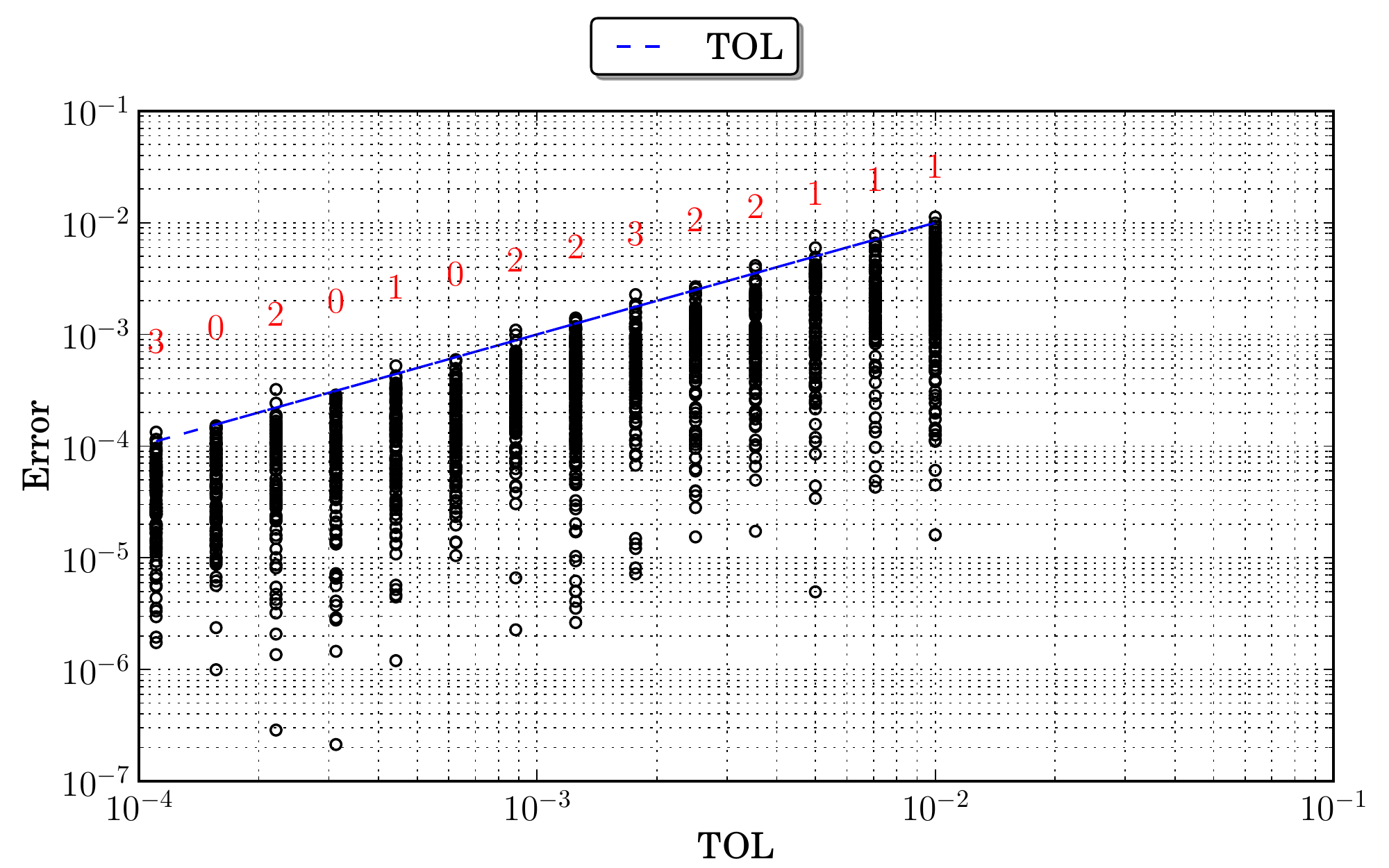}
\end{center}
\caption{From top: {\bf Ex.1, Ex.2, Ex.3}.
    Actual computational errors based on the reference solutions when using CMLMC.
        The numbers above the dashed line show the percentage of runs that
        had errors larger than the required tolerance. We observe that
        \mychanges{these results are consistent with the imposed error
          constraints with a 95\% confidence}.}
\label{fig:ex-error-geo}
\end{figure}

\begin{figure}[ht]
\begin{center}
\includegraphics[scale=0.45]{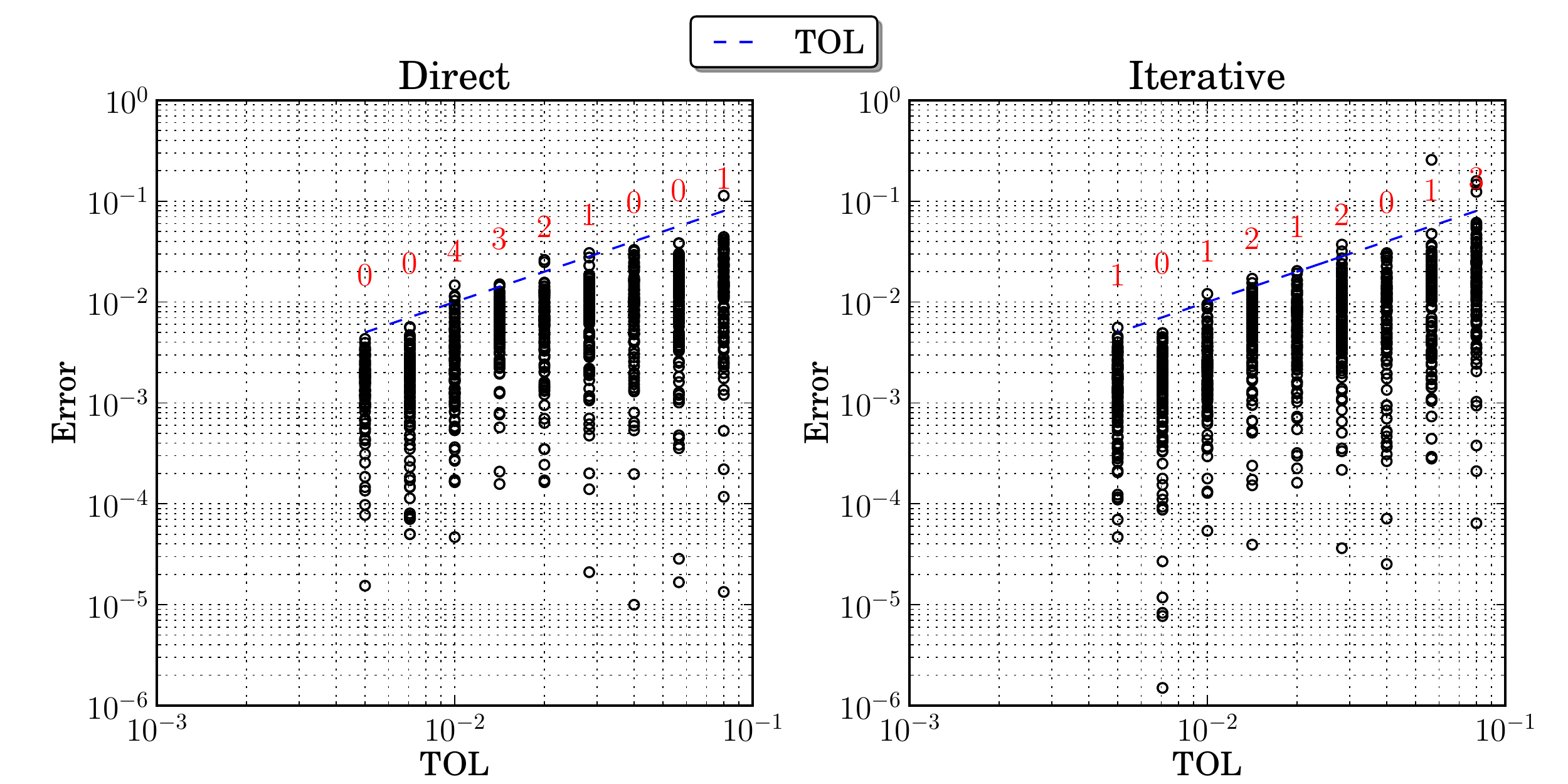}
\includegraphics[scale=0.45]{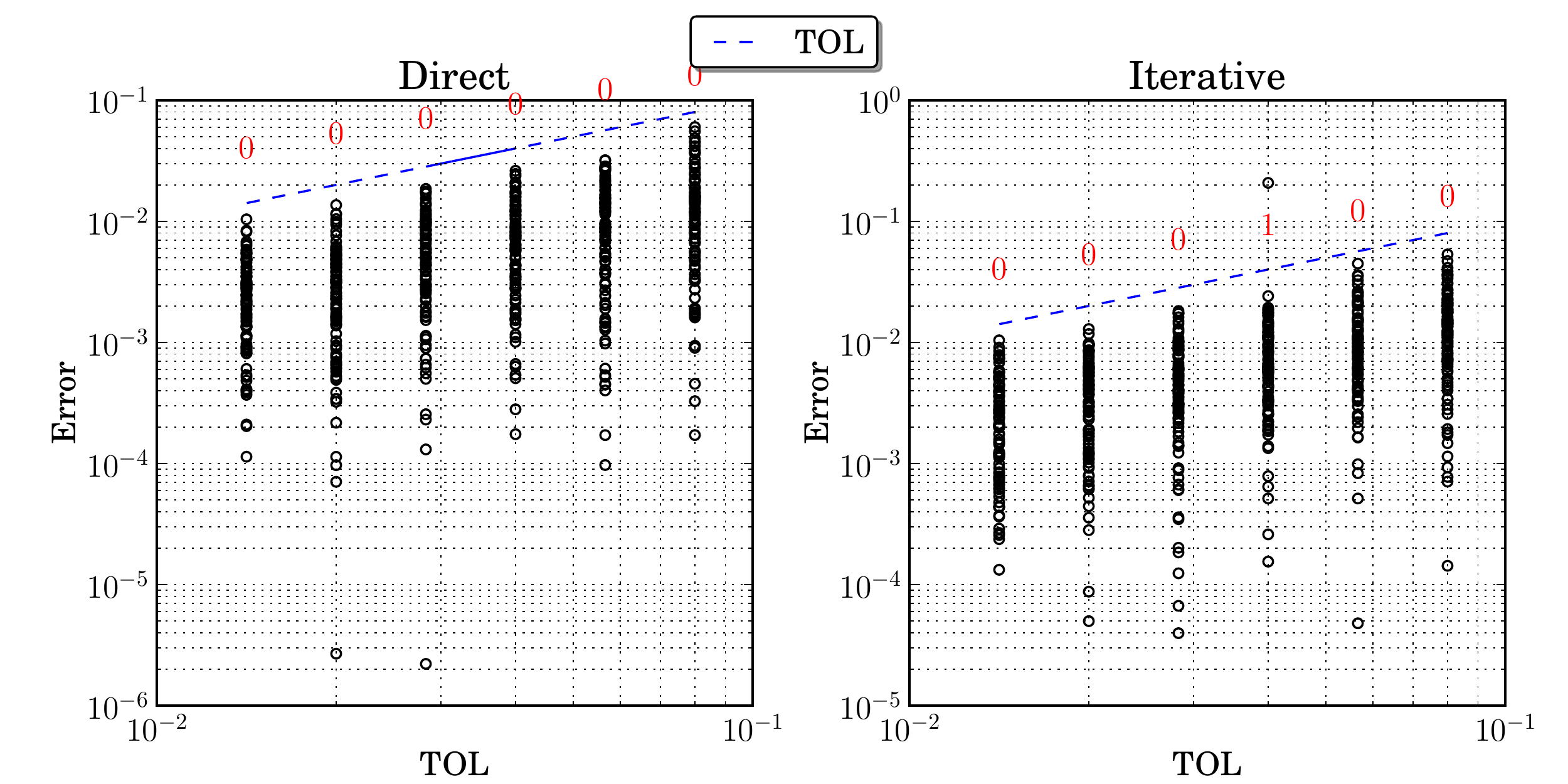}
\includegraphics[scale=0.35]{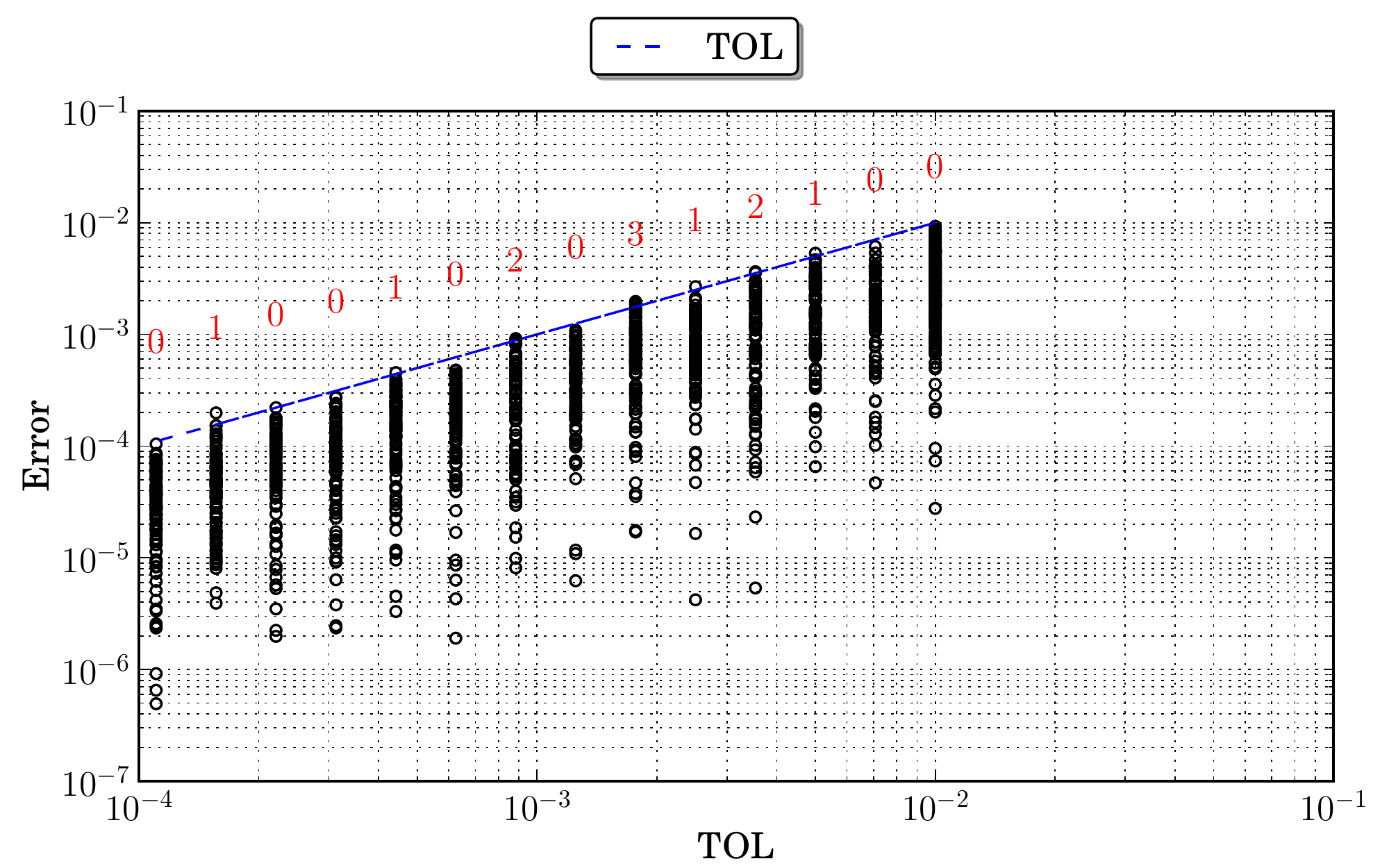}
\end{center}
\caption{From top: {\bf Ex.1, Ex.2, Ex.3}.
        Actual computational errors based on the reference solutions when using SMLMC.
        The numbers above the dashed line show the percentage of runs that
        had errors larger than the required tolerance.
       \mychanges{We observe that these results are consistent with the imposed error
          constraints with a 95\% confidence.
        However, for particular tolerances, the error is smaller than
        $\tol$ because the statistical error is not relaxed when the
        bias is small since tolerance splitting parameter, $\theta$, is
        kept fixed for all tolerances.}}
\label{fig:ex-error-std}
\end{figure}

\begin{figure}[ht]
\begin{center}
\includegraphics[scale=0.45]{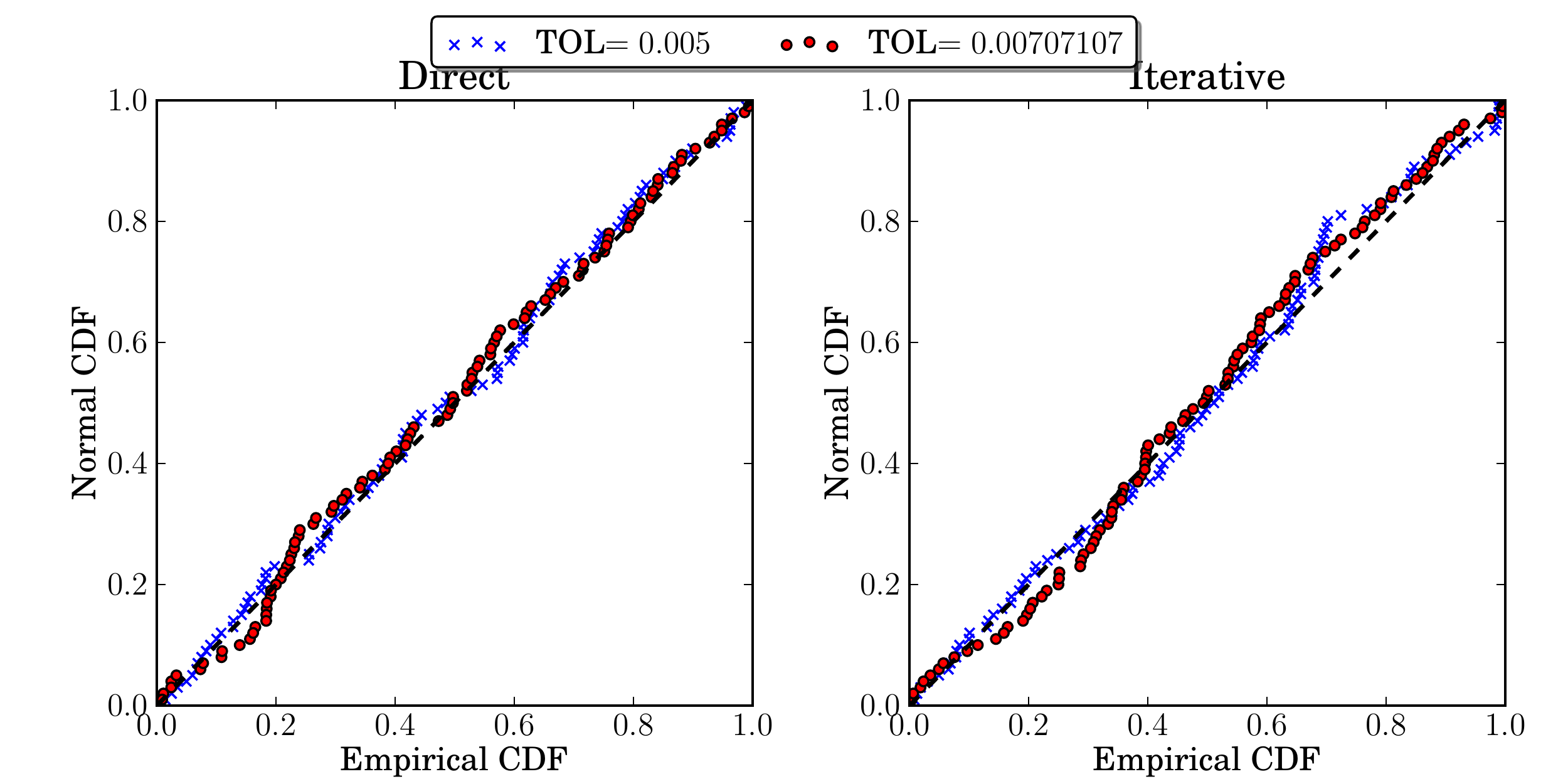}
\includegraphics[scale=0.45]{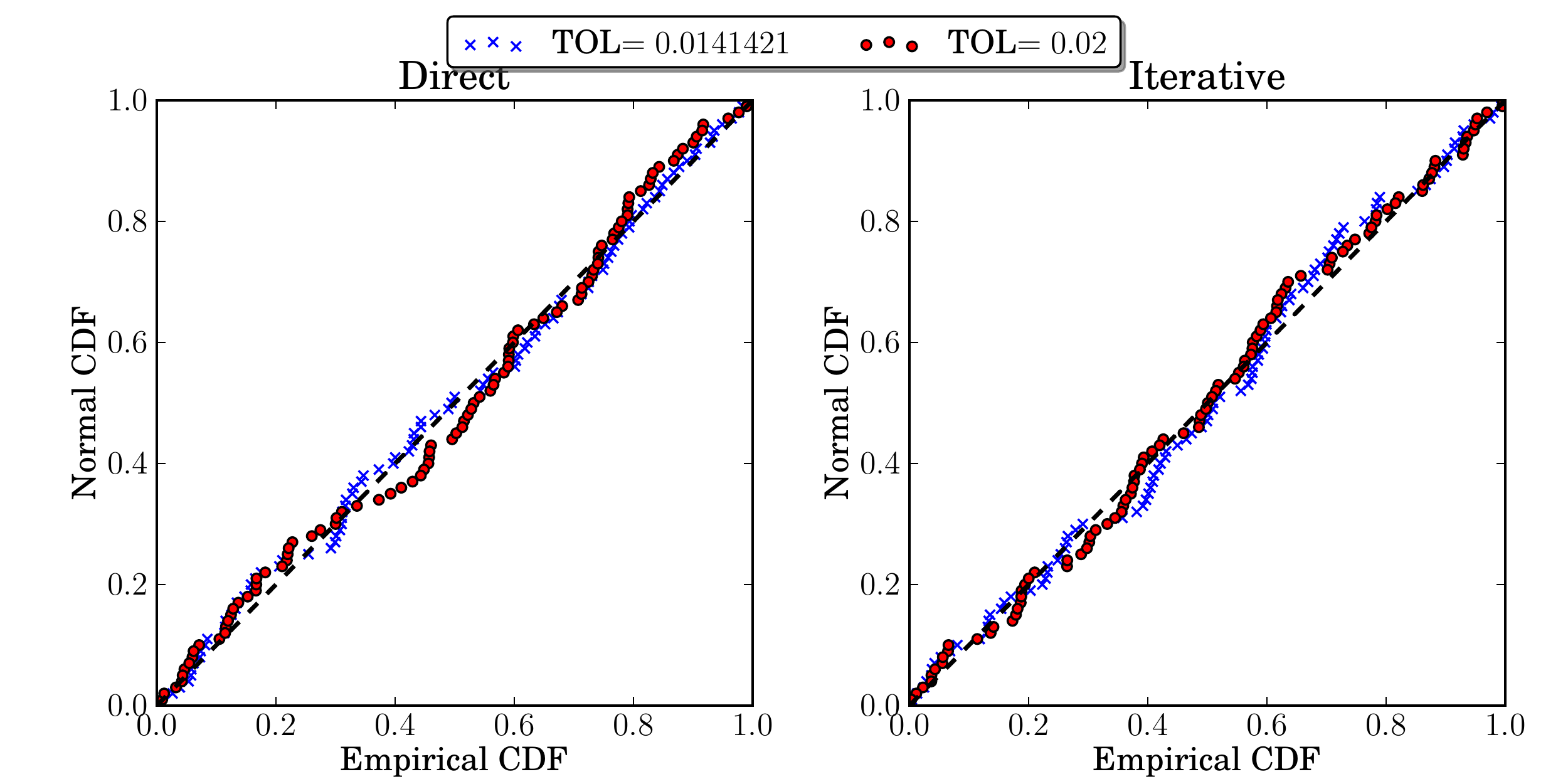}
\includegraphics[scale=0.33]{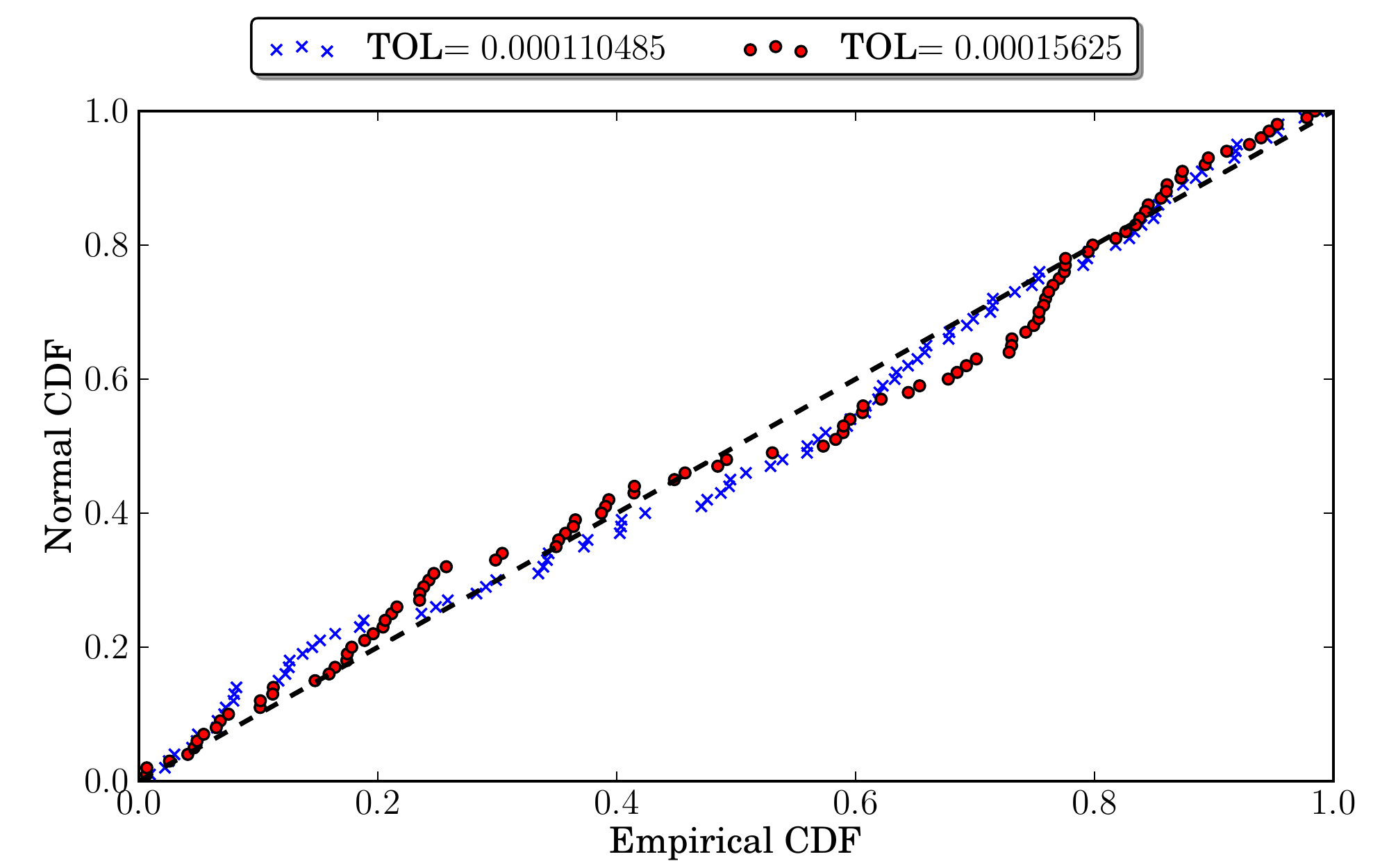}
\end{center}
\caption{From top: {\bf Ex.1, Ex.2, Ex.3}. Normalized empirical
  cumulative distribution function (CDF) of MLMC error estimates for
  different tolerances versus the standard normal CDF. Notice that,
  even for finite tolerances, the standard normal CDF is a good
  approximation of the CDF of the error estimates.}
\label{fig:err_densitites}
\end{figure}

Next, Figure~\ref{fig:ex-error-geo} shows the computational errors of
CMLMC that were computed using the reference solutions as listed in
Section~\ref{ss:examples}. This indicates that the imposed accuracy is
achieved with the required confidence of $95\%$ -- since $C_\alpha=2$.
Compare this figure to Figure~\ref{fig:ex-error-std} which shows the
computational errors of SMLMC. One can see that, in certain cases,
SMLMC solves the problem for a smaller tolerance than the imposed
$\tol$. This is because $\theta$ is fixed and the statistical error is
not relaxed when the bias is small. This can be especially seen in
\textbf{Ex.2} where the choice $h_0=1/8$ produces a bias much smaller
than $0.5\tol$ for the shown tolerances. On the other hand, Figure
\ref{fig:err_densitites} is a QQ-plot showing that the empirical
cumulative distribution function (CDF) of the MLMC error estimates is
well approximated by the standard normal CDF, even for finite
tolerances.

\begin{figure}[ht]
\begin{center}
\includegraphics[scale=0.44]{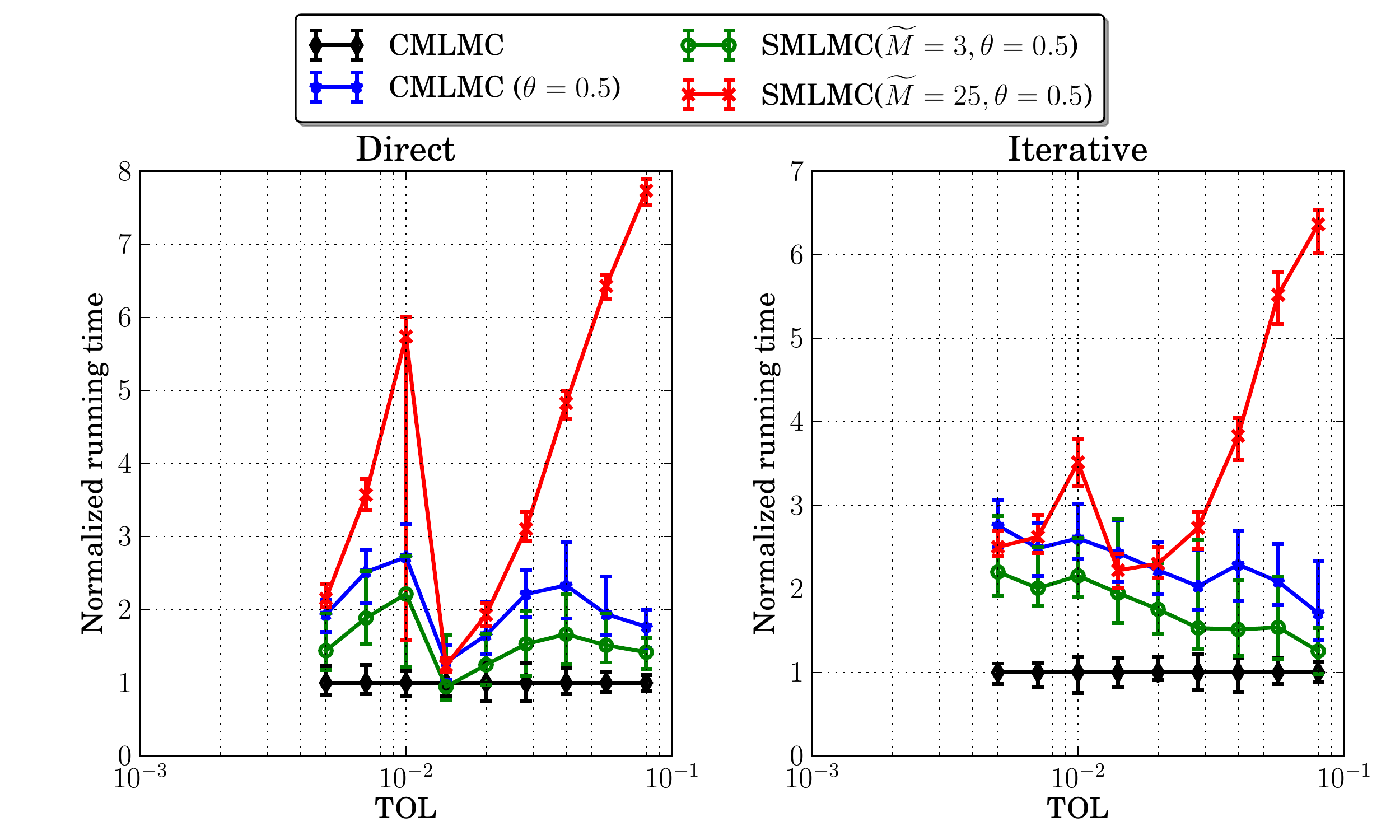}
\includegraphics[scale=0.44]{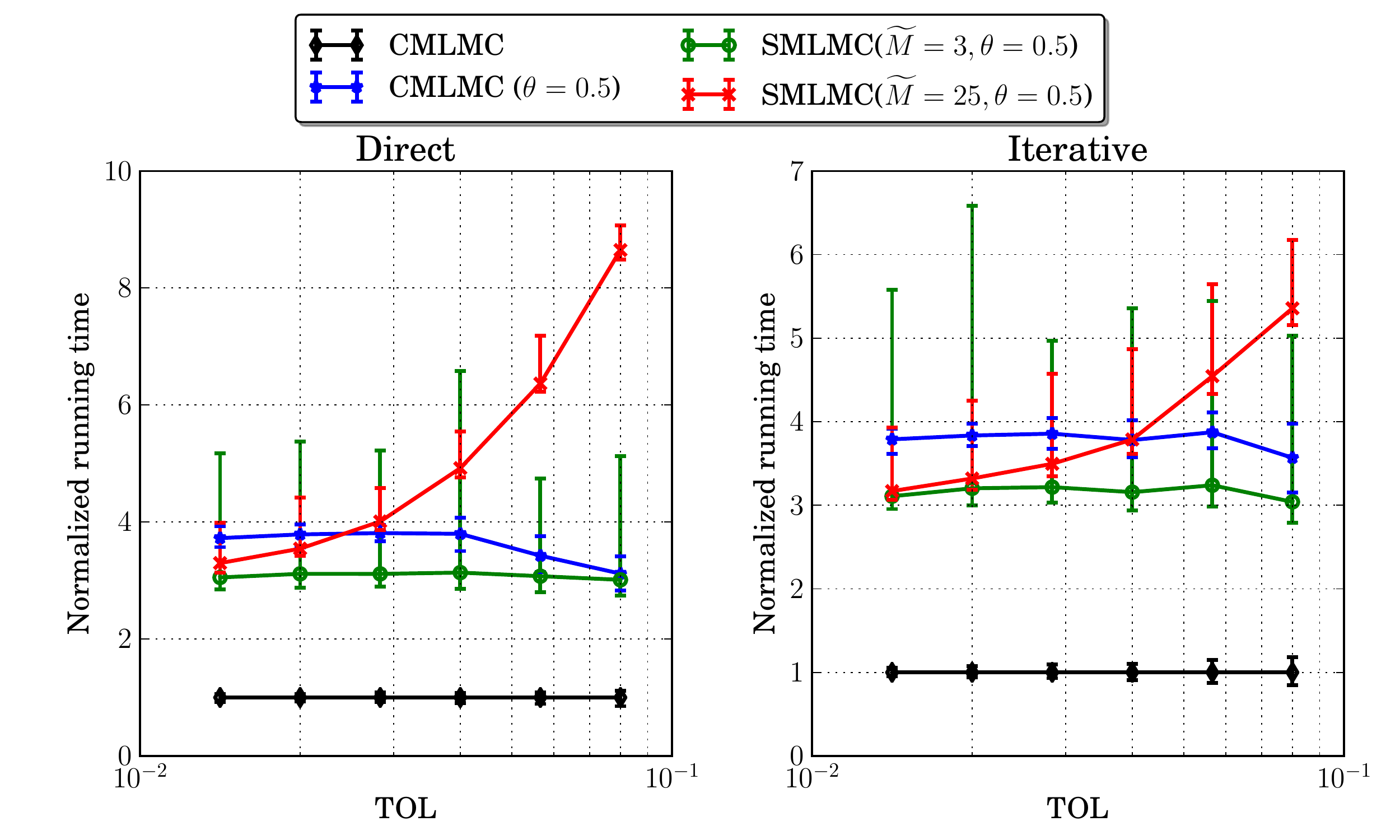}
\includegraphics[scale=0.35]{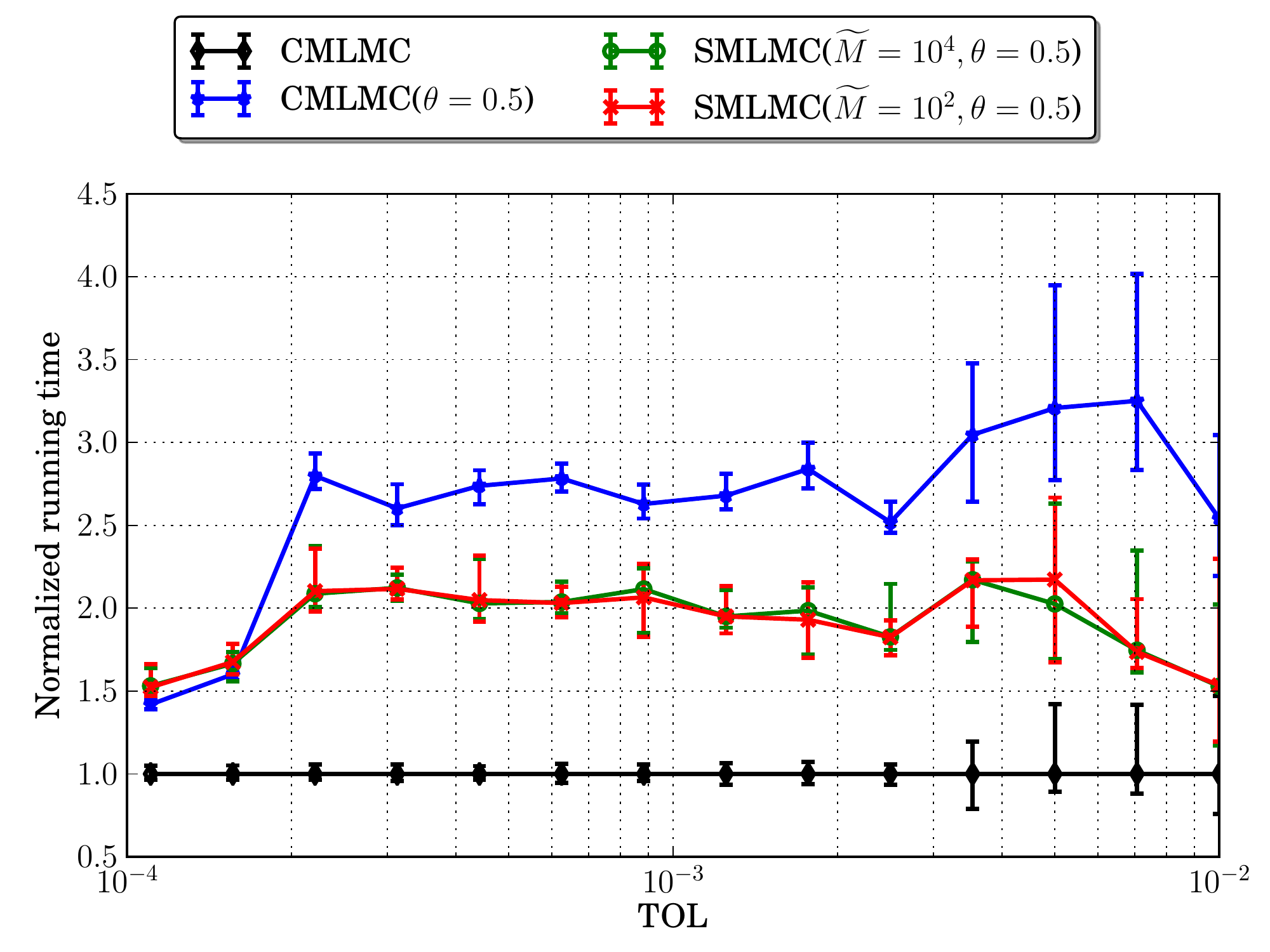}
\end{center}
\caption{From top: {\bf Ex.1, Ex.2, Ex.3}. The running time of CMLMC
  and SMLMC for different $\widetilde M$ and $\theta$, normalized by
  the median running time of CMLMC. This plot shows that a larger
  $\widetilde M$ increases the median running time of the SMLMC but
  also decreases its variability. One sees that CMLMC outperforms
  SMLMC even for a small $\widetilde M$ in all numerical examples.
  \mychanges{Note that for {\bf Ex.1} using direct method and for
    $\tol \approx 0.015$, all algorithms perform similarly. This is
    because, for this example, the optimal error splitting parameter,
    $\theta$, according to \eqref{eq:optimal_theta} is approximately
    0.5.}}
\label{fig:runtime-const-alg}
\end{figure}

\begin{figure}[ht]
\begin{center}
\includegraphics[scale=0.45]{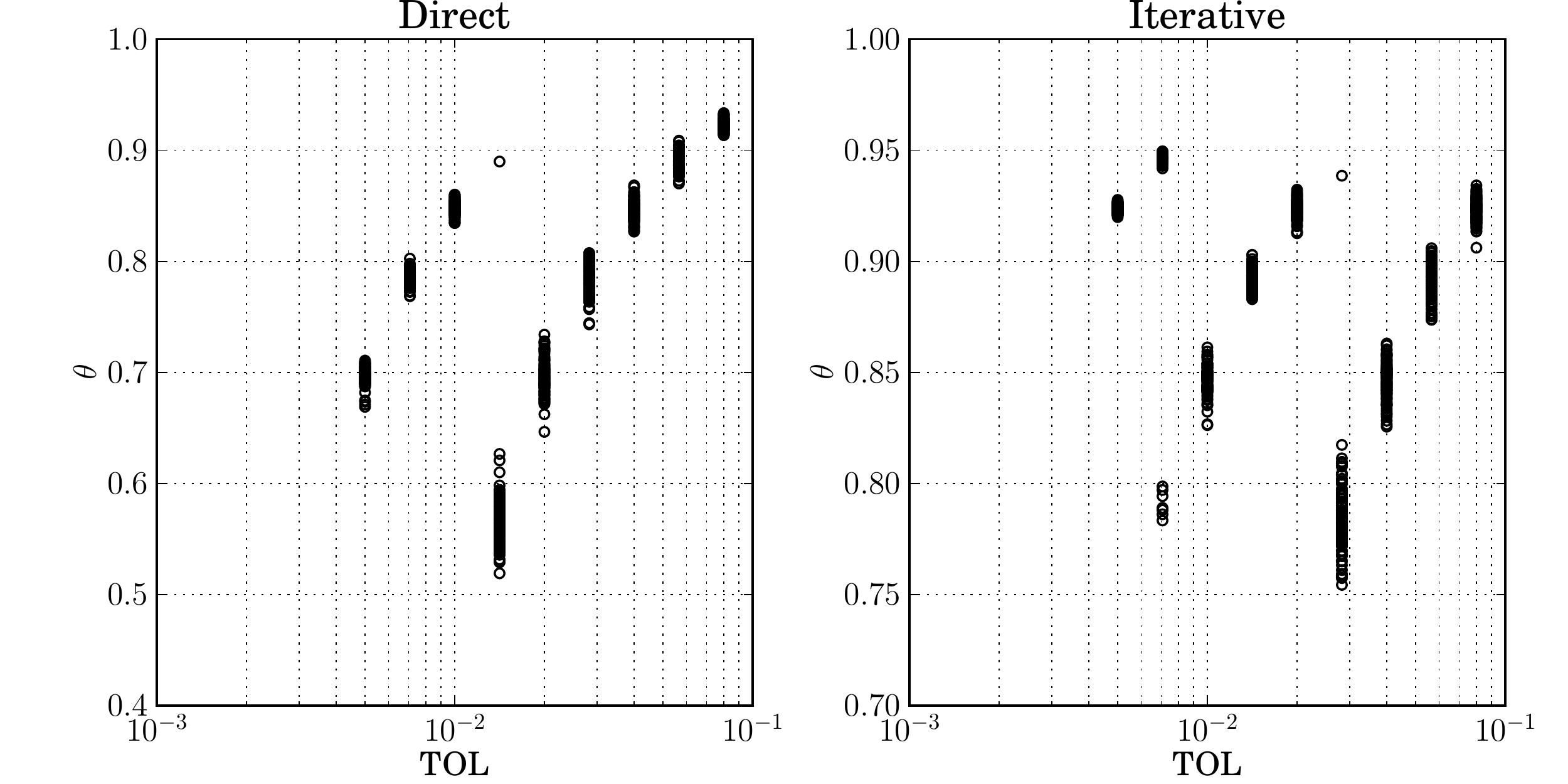}
\includegraphics[scale=0.45]{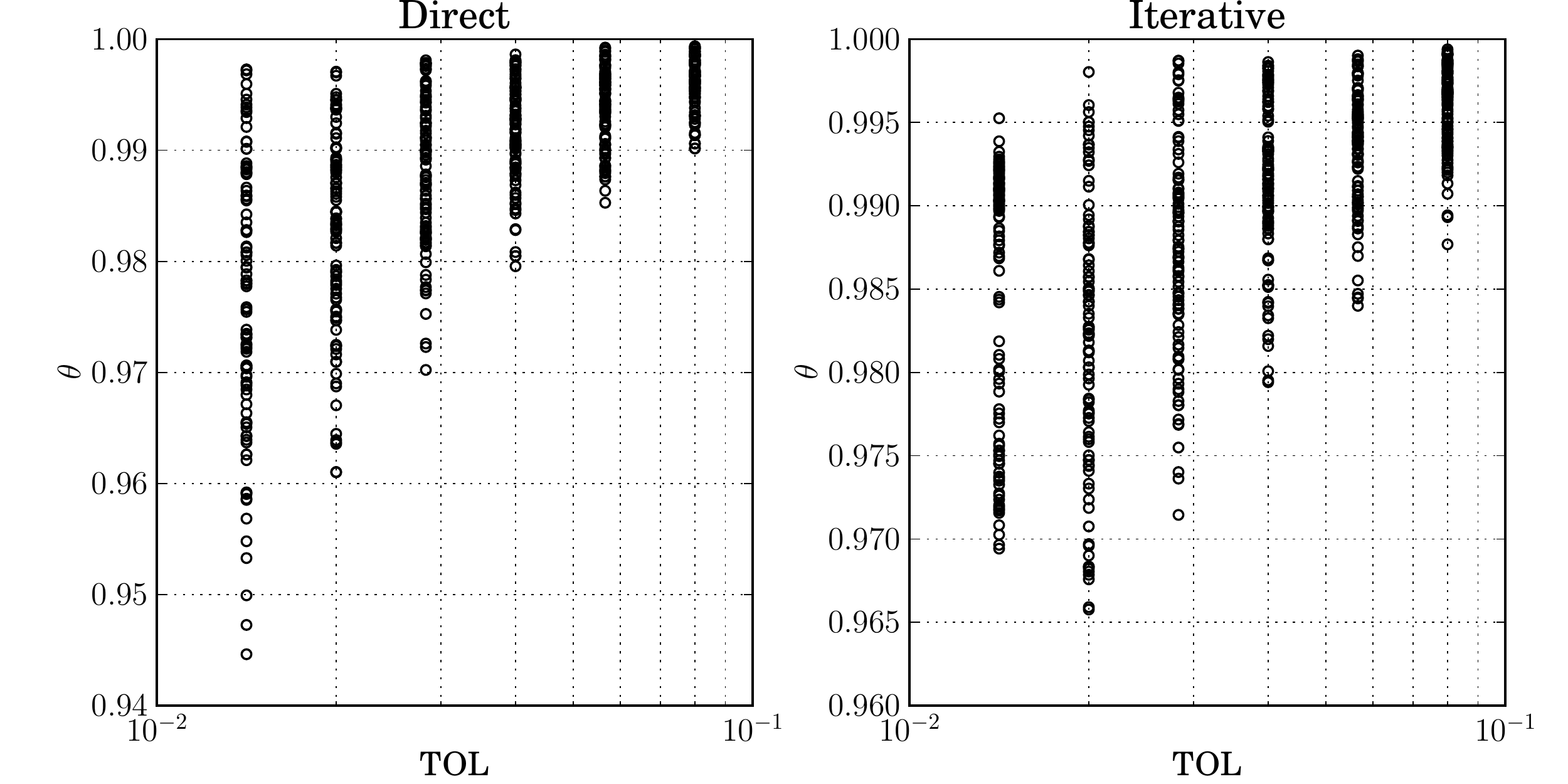}
\includegraphics[scale=0.45]{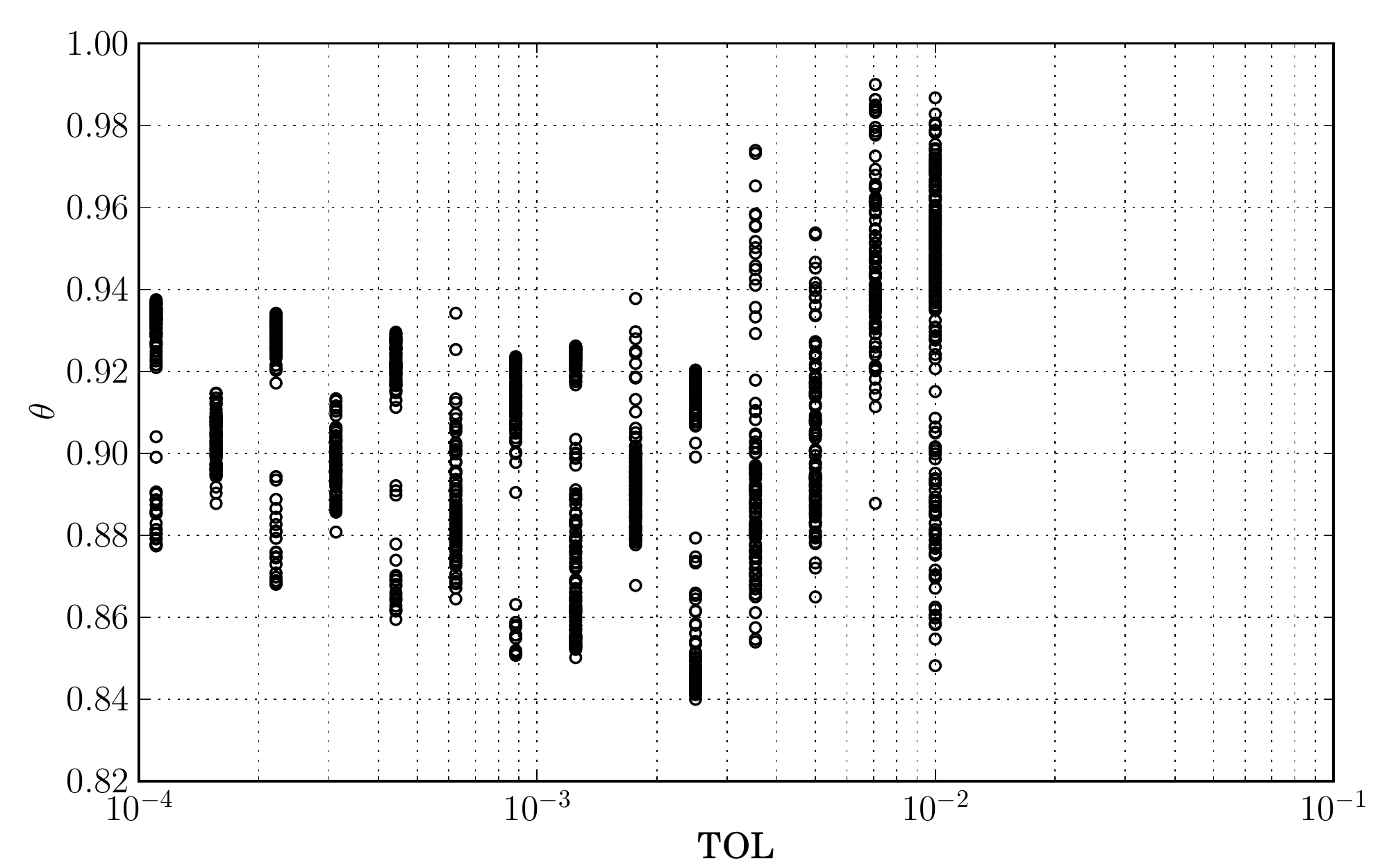}
\end{center}
\caption{From top: {\bf Ex.1, Ex.2, Ex.3}. The error splitting, $\theta$, as computed in \eqref{eq:optimal_theta} an used in CMLMC,
    versus $\tol$.  Observe the behavior of $\theta$ is non-trivial
    and can be far from $\frac{1}{2}$.}
\label{fig:theta-vs-tol}
\end{figure}

Figure~\ref{fig:runtime-const-alg} shows a comparison of the running
time of CMLMC and SMLMC. Notice that a good value of $\widetilde M$ in
SMLMC is not known a priori and the computational time varies
considerably for different values of $\widetilde M$, especially for
smaller tolerances in \textbf{Ex.1} and \textbf{Ex.2}. Specifically, a
larger $\widetilde M$ in SMLMC increases the computational time of the
algorithm, but also decreases its variability. A smaller $\widetilde
M$ gives a smaller computational time at the expense of increased
variation. The variation of the running time is due to inaccurate
estimates of $V_\ell$ due to the smaller number of initial samples. On
the other hand, the running time of CMLMC is less varied, which is a
reflection of the stability of the estimates of $V_\ell$. The
computational savings of CMLMC over SMLMC is an aggregate effect of
the different improvements. This includes 1) a more stable variance
and bias estimates as already discussed, 2) a better splitting of bias
and statistical tolerances. This second point can be seen in
Figure~\ref{fig:theta-vs-tol}, which shows the tolerance splitting
parameter, $\theta$, used in CMLMC as computed by
\eqref{eq:optimal_theta}. We can clearly see here that $\theta$ is not
trivial and changes with the tolerance. Looking closely, one can
notice sudden jumps in the values of $\theta$ due to changes in the
discrete number of levels, $L$. Between jumps, $\theta$ changes
continuously due to inaccuracies in the estimation of the weak error
constant, $Q_W$. Specifically, notice that for $\tol \approx 0.015$ in
\textbf{Ex.1} when using the direct solver, the splitting parameter
$\theta$ used in CMLMC is very close to $0.5$ which explains why, for
this case, the computational time of SMLMC is very close to the
computational time of CMLMC as shown in
Figure~\ref{fig:runtime-const-alg}.

\begin{figure}[ht]
\begin{center}
\includegraphics[scale=0.45]{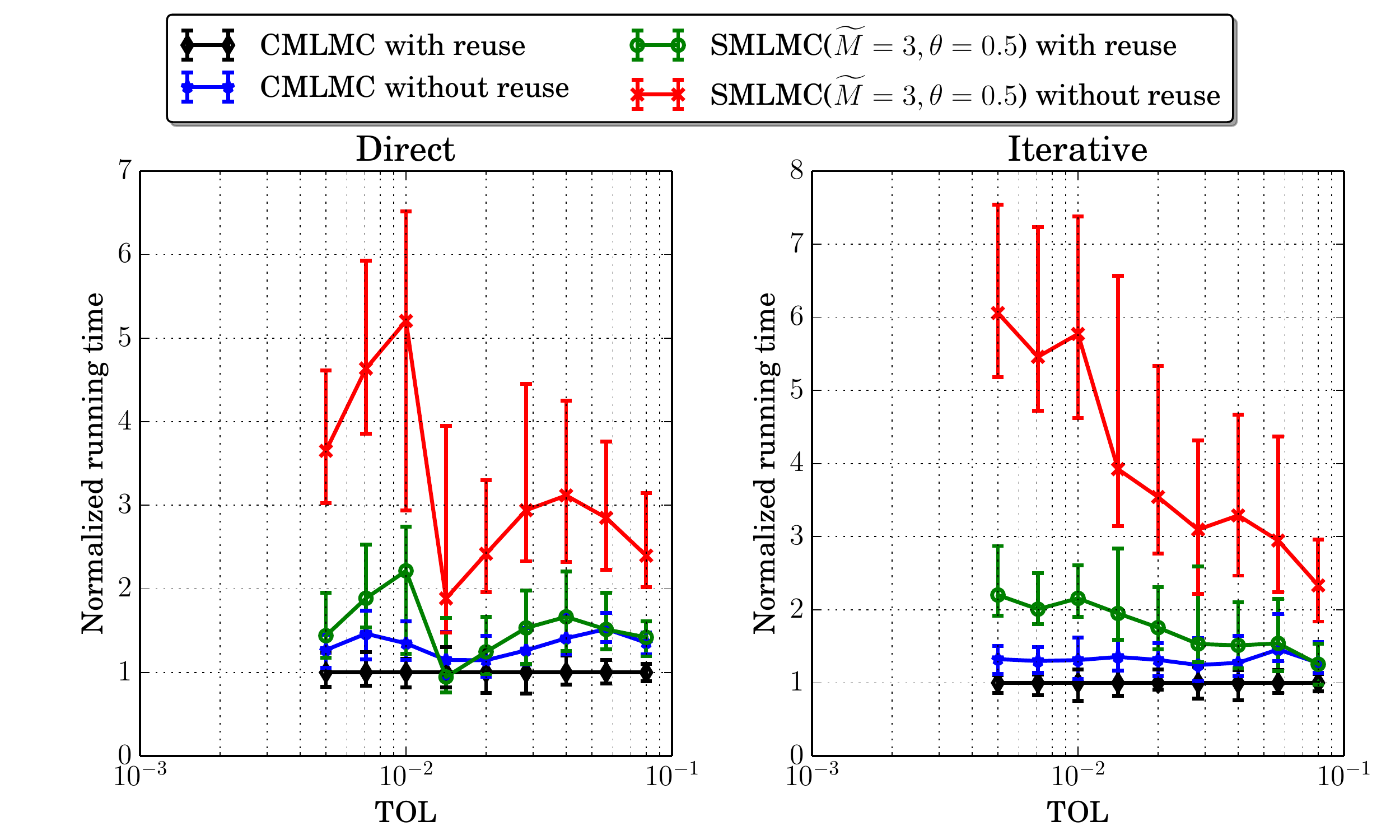}
\includegraphics[scale=0.45]{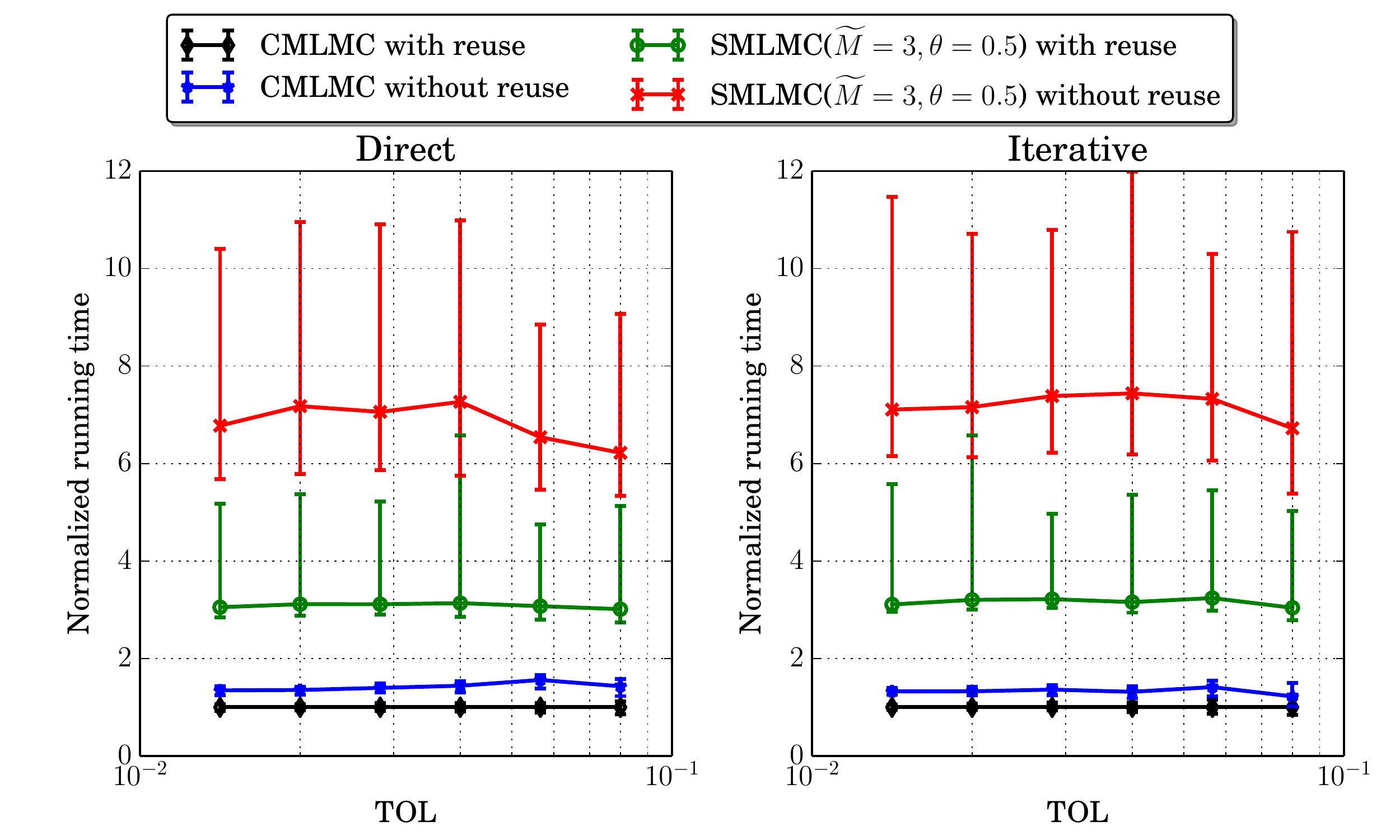}
\includegraphics[scale=0.33]{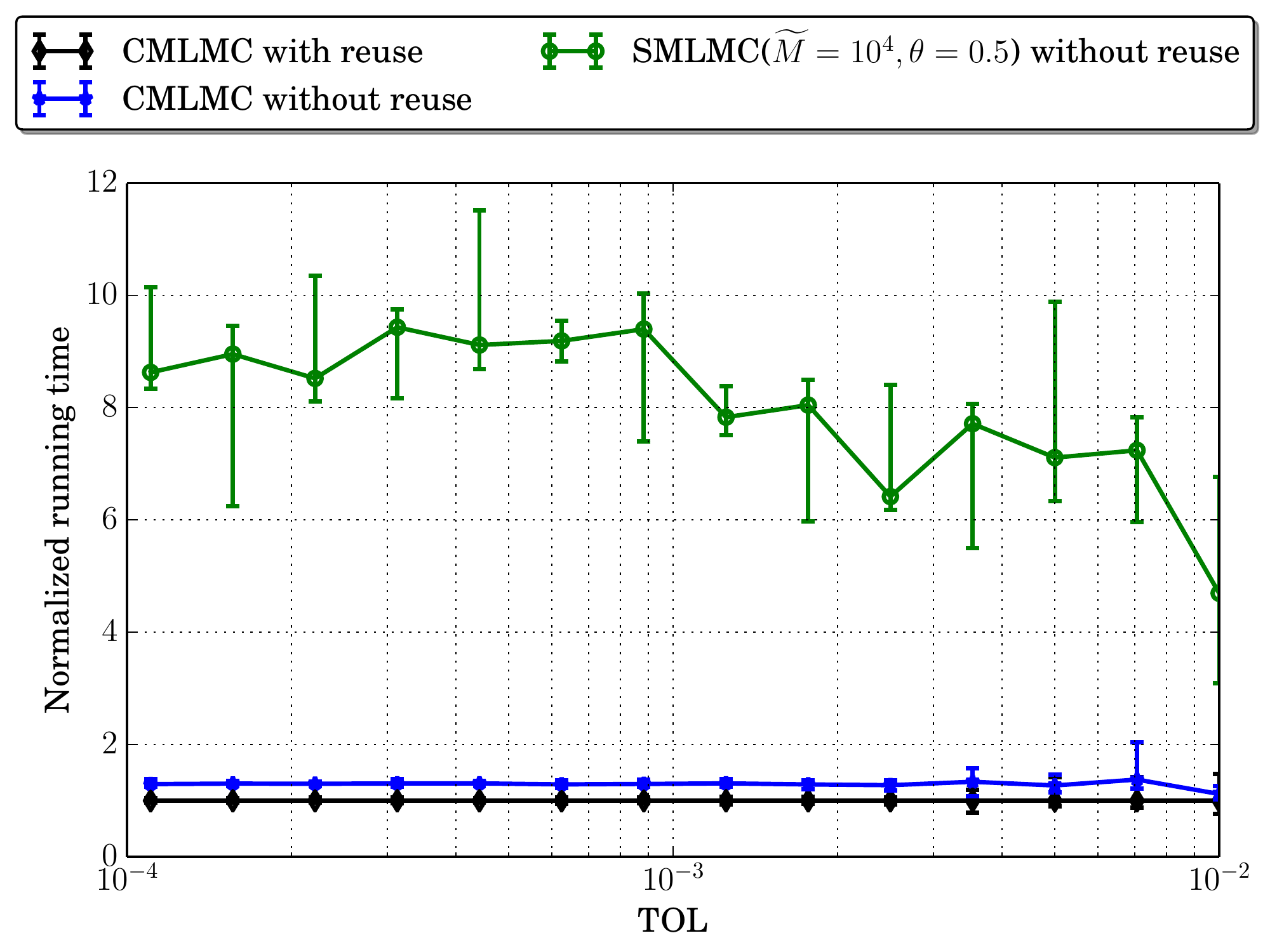}
\end{center}
\caption{From top: {\bf Ex.1, Ex.2, Ex.3}.
    Running time of CMLMC versus SMLMC \mychanges{when} reusing samples for both.
            Also included, is CMLMC without reusing samples from
            previous iterations. All running times are normalized by
            the median of the running time of CMLMC when reusing samples.
            Notice that reusing samples in CMLMC does not add a significant advantage.
            Moreover, CMLMC still produces savings over SMLMC, even when reusing samples in the latter.}
\label{fig:runtime-const-reuse}
\end{figure}

    Finally, the bias of the MLMC estimator when using samples
    generated in previous iterations to compute the quantity of
    interest is not well understood. Using CMLMC, generating new
    samples at each iteration, instead of using samples from previous
    iterations, does not add a significant overhead to the total
    running time of the algorithm.
    Figure~\ref{fig:runtime-const-reuse} explains this point by
    comparing the running time of CMLMC for both cases for both CMLMC
    and SMLMC. This figure shows that computational savings of CMLMC
    over SMLMC whether we reuse samples or not in the former, mainly
    due to better splitting of the tolerance between bias and
    statistical errors. Moreover, it shows that reusing samples in
    CMLMC does not offer significant computational savings that
    justify the increased complexity in the analysis of the resulting
    estimator.

\section{Conclusions}\label{s:conc}  % conc.tex
We have proposed a novel Continuation Multi Level Monte Carlo (CMLMC)
algorithm for weak approximation of stochastic models. Our algorithm uses
discretization hierarchies that are defined a priori for each level
and are geometrically refined across levels. These hierarchies are
either uniform at each level or obtained by regular subdivision of a
non-uniform mesh.

The actual choice of computational work across levels uses the optimal
number of samples per level given the variance and the work
contribution from each level. Accurate computation of these relevant
quantities is based on parametric models.

These parameters are calibrated using approximate samples, either
produced before running the CMLMC and/or during the actual runs. We
also propose a novel Bayesian estimation of the strong and weak error
model parameters, taking particular notice of the deepest levels of
the discretization hierarchy, where only a few realizations are
available to produce the required estimates. The idea is to use
results from coarser levels, where more samples are available, to
stabilize the estimates in the deeper levels. The resulting MLMC
estimator exhibits a non-trivial splitting between bias and
statistical contributions. Indeed, the actual split depends on the
given accuracy and other problem parameters. In fact, as the numerical
examples show, there are cases where most of the accuracy budget is
devoted to the statistical error. Finally, using the Lindeberg-Feller
theorem, we also show the asymptotic normality of the statistical
error in the MLMC estimator and justify in this way our error estimate
that allows prescribing both required accuracy and confidence in the
final result.

We presented three numerical examples to substantiate the above
results, exhibiting the robustness of the new CMLMC Algorithm and to
demonstrate its corresponding computational savings. The examples are
described in terms of differential equations either driven by random
measures or with random coefficients.

Other aspects of MLMC estimators can also be explored, such as the
optimality of geometric hierarchies compared to non-geometric ones.
This will be the subject of a forthcoming work, where extensions of
the CMLMC to that setting will be considered.

\appendix
 % app.tex

\section{Normality of MLMC estimator}
\begin{theorem}\cite[Lindeberg-Feller Theorem, p. 114]{Durret1996}
    For each $n$, let $X_{n,m}$, for $1 \leq n \leq m$, be independent random variables (not necessarily identical).
    Denote
    \begin{align*}
        a_n &= \sum_{m=1}^n X_{n,m}, \\
        Y_{n,m} &=X_{n,m} - \E{X_{n,m}}, \\
        s_n^2 &= \sum_{m=1}^n \E{Y_{n,m}^2}.
    \end{align*}
    Suppose the following Lindeberg condition is satisfied for all $\epsilon>0$:
         \begin{equation}
            \lim_{n \to \infty} s_n^{-2}\sum_{m=1}^n \E{Y_{n,m}^2 \mathbf{1}_{|Y_{n,m} | > \epsilon s_n}} = 0.
            \label{eq:lindeberg_cond}
        \end{equation}
    Then,
        \[ \lim_{n \to \infty} \prob{\frac{a_n - \E{a_n}}{s_n} \leq z } = \Phi(z), \]
    where $\Phi (z)$ is the normal cumulative density function of a standard normal random variable.
\end{theorem}

\begin{lemma}
  \label{thm:clt_result}
  \begin{changes}
  Consider the MLMC estimator $\mathcal{A}$ given by
    \begin{equation*} \mathcal{A} = \sum_{\ell=0}^L \sum_{m=1}^{M_\ell} \frac{G_{\ell}(\omega_{\ell, m})}{M_\ell}, \end{equation*}
    where $G_{\ell}(\omega_{\ell, m})$ denote as usual i.i.d. samples of the random variable $G_\ell$. The family  of random variables, $(G_\ell)_{\ell\ge 0}$, is also assumed independent.
      Denote ${Y_\ell = |G_\ell - \E{G_\ell}|}$ and assume the following
    \begin{subequations}
    \begin{align}
        \label{eq:2nd-moment-bound} C_1 \beta^{-q_3\ell} &\leq
        \E{Y_\ell^2}            & \text{ for all } \ell \geq 0, \\
        \label{eq:4th-moment-bound} \E{Y_\ell^{2+\delta}} &\leq C_2
        \beta^{-\tau \ell} & \text{ for all } \ell \geq 0,
    \end{align}
    \label{eq:cond}
    \end{subequations}
    for some $\beta > 1$ and strictly positive constants $C_1, C_2, q_3,\delta$ and $\tau$.
  Choose the number of samples on each level $M_\ell$ to satisfy,
  for $q_2 > 0$ and a strictly positive sequence $\{H_\ell\}_{\ell \geq 0}$
  \begin{align}
    M_\ell \geq
     \beta^{-q_2\ell} \tol^{-2}  H_\ell^{-1} \left( \sum_{\ell=0}^L
       H_\ell \right) & \qquad \text{for all } \ell \geq 0.
 \label{eq:chosen_ml}
\end{align}
Moreover, choose the number of levels $L$ to satisfy
  \begin{align} L  &\leq \max\left(0, \frac{c \log\left(\tol^{-1}\right)}{\log{\beta}} + C\right) \label{eq:L_cond}
  \end{align}
  for some constants $C$, and $c>0$.
  Finally, denoting
\[ p= (1+\delta/2) q_3 + (\delta/2)q_2 -\tau,\]
  if we have that either $p>0$ or $c < \delta/p$,
  then
  \[\lim_{\tol \to 0} \prob{\frac{\mathcal{A}- \E{\mathcal
        A}}{\sqrt{\var{\mathcal A}}} \leq z} = \Phi \left(z
  \right). \]
  \end{changes}
\end{lemma}
\proof
    We prove this lemma by ensuring that the Lindeberg condition \eqref{eq:lindeberg_cond} is satisfied.
    The condition becomes in this case
    \[ \lim_{\tol \to 0} \underbrace{\frac{1}{\var{\mathcal A}} \sum_{\ell=0}^L \sum_{m=1}^{M_\ell} \E{ \frac{Y_\ell^2}{M_\ell^2} \mathbf{1}_{\frac{Y_\ell}{M_\ell}
        > \epsilon \sqrt{\var{\mathcal A}}}}}_{:= F} = 0,\]
    for all $\epsilon > 0$. Below we make repeated use of the following identity for non-negative sequences ${\{a_\ell \}}$ and ${\{b_\ell \}}$ and $q \geq 0$.
    \begin{equation} \sum_{\ell}{a_\ell^q b_\ell} \leq \left( \sum_\ell a_\ell\right)^q \sum_\ell b_\ell. \label{eq:identity} \end{equation}
    First we use the Markov inequality to bound
    \begin{align*}
        F &=  \frac{1}{\var{\mathcal A}} \sum_{\ell=0}^L \sum_{m=1}^{M_\ell} \E{ \frac{Y_\ell^2}{M_\ell^2} \mathbf{1}_{Y_\ell >
                    \epsilon \sqrt{\var{\mathcal A}} M_\ell}} \\
        &\leq \frac{\epsilon^{-\delta}}{\var{\mathcal A}^{1+\delta/2}}  \sum_{\ell=0}^L  M_\ell^{-1-\delta} \E{Y_\ell^{2+\delta}}.
    \end{align*}
    Using \eqref{eq:identity} and substituting for the variance
    $\var{\mathcal A}$ where we denote $\var{G_\ell} = \E{\left(
        G_\ell - \E{G_\ell} \right)^{2}}$ by $V_\ell$, we find
    \begin{changes}
    \begin{align*}
        F &\leq \frac{\epsilon^{-\delta} \left( \sum_{\ell=0}^L  M_\ell^{-1} V_\ell \right)^{1+\delta/2}}{\left( \sum_{\ell=0}^L V_\ell M_\ell^{-1} \right)^{1+\delta/2}}
                \sum_{\ell=0}^L V_\ell^{-1-\delta/2} M_\ell^{-\delta/2} \E{Y_\ell^{2+\delta}} \\
        &\leq \epsilon^{-\delta}
                \sum_{\ell=0}^L V_\ell^{-1-{\delta}/{2}} M_\ell^{-{\delta}/{2}} \E{Y_\ell^{2+\delta}}.
    \end{align*}
     Using the lower bound on the number of samples $M_\ell$ \eqref{eq:chosen_ml} and \eqref{eq:identity} again yields
    \begin{align*}
        F &\leq \epsilon^{-\delta} \tol^\delta \left(   \sum_{\ell=0}^L V_\ell^{-1-{\delta}/{2}} \beta^{\frac{\delta q_2 \ell}{2}} H_\ell^{\delta/2} \E{Y_\ell^{2+\delta}} \right)
                \left( \sum_{\ell=0}^L H_\ell \right)^{-\delta/2} \\
        &\leq \epsilon^{-\delta} \tol^\delta \left( \sum_{\ell=0}^L V_\ell^{-1-{\delta}/{2}} \beta^{(\delta/2)q_2 \ell} \E{Y_\ell^{2+\delta}} \right).
    \end{align*}
    Finally using the bounds \eqref{eq:2nd-moment-bound} and \eqref{eq:4th-moment-bound}
    \begin{align*}
      F &\leq \epsilon^{-\delta} \tol^\delta \left(C_1^{-1-{\delta}/{2}} C_2 \sum_{\ell=0}^L \beta^{(1+\delta/2) q_3\ell} \beta^{(\delta/2)q_2 \ell} \beta^{-\tau \ell} \right)\\
      &= \epsilon^{-\delta}\tol^{\delta} C_1^{-1-{\delta}/{2}} C_2
      \frac{\beta^{(L+1)p} - 1}{\beta^p-1} ,
    \end{align*}
    We distinguish two cases here, namely:
    \begin{itemize}
    \item If $p>0$ is satisfied then $\lim_{\tol \to 0} F= 0$ for any $c>0$.
    \item Otherwise, substituting \eqref{eq:L_cond} gives
    \begin{align*}
        F &\leq \epsilon^{-\delta} \tol^{\delta}
        C_1^{-1-{\delta}/{2}} C_2  \frac{\tol^{-cp } \beta^{(C+1)p} -
          1}{\beta^p-1} = \Order{\tol^{\delta - cp}},
    \end{align*}
    and since in this case $cp < \delta$ then $\lim_{\tol \to 0} F = 0$. \hfill $\square$
    \end{itemize}
\end{changes}

\begin{remark}
    The choice \eqref{eq:chosen_ml} mirrors the choice \eqref{eq:optimal_ml}, the latter being the optimal number of samples to bound the statistical error of the estimator by $\tol$.
    Specifically, $H_\ell \propto \sqrt{V_\ell W_\ell}$ where $W_\ell$ is the work per sample on level $\ell$.
    Moreover, the choice \eqref{eq:optimal_ml} uses the variances $\{ V_\ell \}_{\ell=0}^L$ or an estimate of it in the actual implementation.
    On the other hand, the choice \eqref{eq:chosen_ml} uses the upper bound of $V_\ell$ instead, if $q_2$ is the rate of strong convergence therein.
    Furthermore, if we assume the weak error model \eqref{eq:weak_error_model} holds and $h_L = h_0 \beta^{-L}$
    then we must have
    \[ Q_W h_L^{q_1} = Q_W h_0^{q_1} \beta^{-L q_1} \leq (1-\theta)\tol, \]
    which gives a lower bound on the number of levels $L$, namely
    \[  L \geq \frac{\log(\tol^{-1})}{q_1\log(\beta)} + \frac{-\log(1-\theta) +\log(Q_W) + q_1\log(h_0)}{q_1\log(\beta)}, \]
    to bound the bias by $\tol$.

    Finally, in Example \ref{ex:spde_problem} the conditions
    \eqref{eq:cond} are satisfied for $q_3=2$ and, assuming $p^*> 3$, for $\delta=1$ and $\tau=6$.
    Similarly, Example \ref{ex:sde_problem} satisfies the conditions \eqref{eq:cond} are for $q_3=1$ and $\delta=2$ and $\tau=2$, cf. \cite{hsst12}.
\end{remark}

\begin{remark}
    The assumption \eqref{eq:2nd-moment-bound} can be relaxed. For instance, one can assume instead
    that
    \begin{align*}
        V_{\ell+1} &\leq V_{\ell} \qquad \text{  for all }\ell \geq 1,\\
       0 & < \lim_{\ell \to \infty} \var{Y_\ell} \beta^{q_3 \ell} < \infty,
    \end{align*}
    and slightly different conditions on $L$.
\end{remark}

\begin{acknowledgements} Ra\'{u}l~Tempone is a member  of the
  Strategic Research Initiative on Uncertainty
Quantification in Computational Science and Engineering at KAUST (SRI-UQ).
The authors would like to recognize the support of King Abdullah University of Science and
Technology (KAUST) AEA project ``Predictability and Uncertainty
Quantification for Models of Porous Media'' and University of Texas at Austin AEA Rnd 3
``Uncertainty quantification for predictive modeling of the dissolution of porous and fractured media''.
We would also like to acknowledge the use of the following open source software packages: {\tt PETSc}~\cite{petsc-web-page}, {\tt PetIGA}~\cite{Collier2013}, {\tt NumPy}, {\tt matplotlib}~\cite{matplotlib}.
\end{acknowledgements}

\bibliographystyle{spmpsci}

\begin{thebibliography}{10}
\providecommand{\url}[1]{#1}
\providecommand{\urlprefix}{URL }
\expandafter\ifx\csname urlstyle\endcsname\relax
  \providecommand{\doi}[1]{DOI~\discretionary{}{}{}#1}\else
  \providecommand{\doi}{DOI~\discretionary{}{}{}\begingroup
  \urlstyle{rm}\Url}\fi

\bibitem{Amestoy2001}
Amestoy, P.R., Duff, I.S., L'Excellent, J.Y., Koster, J.: A fully asynchronous
  multifrontal solver using distributed dynamic scheduling.
\newblock SIAM J. Matrix Anal. Appl. \textbf{23}, 15--41 (2001).
\newblock \doi{10.1137/S0895479899358194}.
\newblock \urlprefix\url{http://portal.acm.org/citation.cfm?id=587708.587825}

\bibitem{Amestoy2006}
Amestoy, P.R., Guermouche, A., L'Excellent, J.Y., Pralet, S.: Hybrid scheduling
  for the parallel solution of linear systems.
\newblock Parallel Computing \textbf{32}(2), 136 -- 156 (2006).
\newblock \doi{10.1016/j.parco.2005.07.004}.
\newblock
  \urlprefix\url{http://www.sciencedirect.com/science/article/pii/S0167819105001328}

\bibitem{bnt2010}
Babu{\v{s}}ka, I., Nobile, F., Tempone, R.: A stochastic collocation method for
  elliptic partial differential equations with random input data.
\newblock SIAM review \textbf{52}(2), 317--355 (2010)

\bibitem{petsc-web-page}
Balay, S., Brown, J., Buschelman, K., Gropp, W.D., Kaushik, D., Knepley, M.G.,
  McInnes, L.C., Smith, B.F., Zhang, H.: {PETSc} {W}eb page (2013).
\newblock Http://www.mcs.anl.gov/petsc

\bibitem{bsz11}
Barth, A., Schwab, C., Zollinger, N.: Multi-level {M}onte {C}arlo finite
  element method for elliptic {PDE}s with stochastic coefficients.
\newblock Numerische Mathematik \textbf{119}(1), 123--161 (2011)

\bibitem{cst13}
Charrier, J., Scheichl, R., Teckentrup, A.: Finite element error analysis of
  elliptic {PDE}s with random coefficients and its application to multilevel
  {M}onte {C}arlo methods.
\newblock SIAM Journal on Numerical Analysis \textbf{51}(1), 322--352 (2013)

\bibitem{cgst11}
Cliffe, K., Giles, M., Scheichl, R., Teckentrup, A.: Multilevel {M}onte {C}arlo
  methods and applications to elliptic {PDE}s with random coefficients.
\newblock Computing and Visualization in Science \textbf{14}(1), 3--15 (2011)

\bibitem{Collier2013}
Collier, N., Dalcin, L., Calo, V.: {PetIGA}: High-performance isogeometric
  analysis.
\newblock arxiv (1305.4452) (2013).
\newblock Http://arxiv.org/abs/1305.4452

\bibitem{petiga}
Dalcin, L., Collier, N.: Pet{IGA}: A framework for high performance
  isogeometric analysis (2013).
\newblock Https://bitbucket.org/dalcinl/petiga

\bibitem{Durret1996}
Durrett, R.: Probability: theory and examples.
\newblock second edn. Duxbury Press, Belmont, CA (1996)

\bibitem{giles08_m}
Giles, M.: Improved multilevel {M}onte {C}arlo convergence using the milstein
  scheme.
\newblock In: {M}onte {C}arlo and quasi-{M}onte {C}arlo methods 2006, pp.
  343--358. Springer (2008)

\bibitem{giles08}
Giles, M.: Multilevel {M}onte {C}arlo path simulation.
\newblock Operations Research \textbf{56}(3), 607--617 (2008)

\bibitem{gr12}
Giles, M., Reisinger, C.: Stochastic finite differences and multilevel {M}onte
  {C}arlo for a class of {SPDEs} in finance.
\newblock SIAM Journal of Financial Mathematics \textbf{3}(1), 572--592 (2012)

\bibitem{gs13c}
Giles, M., Szpruch, L.: Antithetic multilevel {M}onte {C}arlo estimation for
  multidimensional {SDE}s.
\newblock In: {M}onte {C}arlo and Quasi-{M}onte {C}arlo Methods 2012
  (submitted). Springer (2013)

\bibitem{gs13b}
Giles, M., Szpruch, L.: Multilevel {M}onte {C}arlo methods for applications in
  finance, pp. 3--48.
\newblock World Scientific (2013)

\bibitem{gs13}
Giles, M., Szpruch, L.: Antithetic multilevel {M}onte {C}arlo estimation for
  multi-dimensional {SDE}s without {L\'e}vy area simulation.
\newblock To appear in Annals of Applied Probability  (2013/4)

\bibitem{Glasserman}
Glasserman, P.: {M}onte {C}arlo methods in financial engineering,
  \emph{Applications of Mathematics (New York)}, vol.~53.
\newblock Springer-Verlag, New York (2004).
\newblock Stochastic Modelling and Applied Probability

\bibitem{heinrich98}
Heinrich, S.: {M}onte {C}arlo complexity of global solution of integral
  equations.
\newblock Journal of Complexity \textbf{14}(2), 151--175 (1998)

\bibitem{hs99}
Heinrich, S., Sindambiwe, E.: {M}onte {C}arlo complexity of parametric
  integration.
\newblock Journal of Complexity \textbf{15}(3), 317--341 (1999)

\bibitem{hsst12}
Hoel, H., Schwerin, E.v., Szepessy, A., Tempone, R.: Adaptive multilevel
  {M}onte {C}arlo simulation.
\newblock In: Engquist, B., Runborg, O., Tsai, Y.H. (eds.) Numerical Analysis
  of Multiscale Computations, no.~82 in Lecture Notes in Computational Science
  and Engineering, pp. 217--234. Springer (2012)

\bibitem{matplotlib}
Hunter, J.D.: Matplotlib: A {2D} graphics environment.
\newblock Computing In Science \& Engineering \textbf{9}(3), 90--95 (2007)

\bibitem{JouCviMus_Book}
Jouini, E., Cvitani{\'c}, J., Musiela, M. (eds.): Option pricing, interest
  rates and risk management.
\newblock Handbooks in Mathematical Finance. Cambridge University Press,
  Cambridge (2001)

\bibitem{KS}
Karatzas, I., Shreve, S.E.: Brownian motion and stochastic calculus,
  \emph{Graduate Texts in Mathematics}, vol. 113.
\newblock Second edn. Springer-Verlag, New York (1991)

\bibitem{kebaier05}
Kebaier, A.: Statistical {R}omberg extrapolation: a new variance reduction
  method and applications to options pricing.
\newblock Annals of Applied Probability \textbf{14}(4), 2681--2705 (2005)

\bibitem{Mordeckietal08}
Mordecki, E., Szepessy, A., Tempone, R., Zouraris, G.E.: Adaptive weak
  approximation of diffusions with jumps.
\newblock SIAM J. Numer. Anal. \textbf{46}(4), 1732--1768 (2008)

\bibitem{Ok}
{\O}ksendal, B.: Stochastic differential equations.
\newblock Universitext, fifth edn. Springer-Verlag, Berlin (1998)

\bibitem{sivia1996data}
Sivia, D.S.: Data Analysis.: A Bayesian Tutorial.
\newblock Oxford University Press (1996)

\bibitem{stz2001}
Szepessy, A., Tempone, R., Zouraris, G.E.: Adaptive weak approximation of
  stochastic differential equations.
\newblock Communications on Pure and Applied Mathematics \textbf{54}(10),
  1169--1214 (2001)

\bibitem{teckentrup13}
Teckentrup, A.: Multilevel {M}onte {C}arlo methods and uncertainty
  quantification.
\newblock {PhD} thesis, University of Bath (2013)

\bibitem{tsgu13}
Teckentrup, A., Scheichl, R., Giles, M., Ullmann, E.: Further analysis of
  multilevel {M}onte {C}arlo methods for elliptic {PDE}s with random
  coefficients.
\newblock Numerische Mathematik \textbf{125}(3), 569--600 (2013)

\bibitem{xg12}
Xia, Y., Giles, M.: Multilevel path simulation for jump-diffusion {SDE}s.
\newblock In: Plaskota, L., Wo{\'z}niakowski, H. (eds.) {M}onte {C}arlo and
  Quasi-{M}onte {C}arlo Methods 2010, pp. 695--708. Springer (2012)

\end{thebibliography}

\end{document}